\documentclass[journal]{IEEEtran}
\usepackage{amssymb}
\ifCLASSINFOpdf
\else
\usepackage[dvips]{graphicx}
\fi
\usepackage{graphics}
\usepackage{subfigure}
\usepackage{graphicx}
\usepackage{algorithm}
\usepackage{algorithmic}
\usepackage{graphics}
\usepackage{subfig}
\usepackage{balance}
\usepackage{amsmath}
\usepackage{color}
\usepackage{colortbl}
\usepackage{xcolor}
\usepackage{multirow}
\newtheorem{theorem}{Theorem}
\newtheorem{definition}{Definition}
\newtheorem{lemma}{Lemma}
\newtheorem{proposition}{Proposition}

\newtheorem{remark}{Remark}
\begin{document}
\title{Inertial nonconvex alternating minimizations for the image deblurring}
\author{Tao Sun, Roberto Barrio, Marcos Rodr\'{\i}guez, and Hao Jiang % <-this % stops a space
\thanks{This work was
supported by the National Key Research and Development Program of China (2017YFB0202003), and National Natural Science Foundation of Hunan Province in China (2018JJ3616), and by the Spanish Research projects MTM2015-64095-P, PGC2018-096026-B-I00, and European Regional Development Fund and Diputaci\'on General de Arag\'on (E24-17R).

Tao Sun and Hao Jiang are with the
College of Computer, National University of Defense Technology,
Changsha, 410073, Hunan,  China (e-mail: \texttt{nudtsuntao@163.com}, \texttt{haojiang@nudt.edu.cn}).

Roberto Barrio is with Departamento de Matem\'atica Aplicada, IUMA and CoDy group, Universidad de Zaragoza, Pedro Cerbuna, 12, Zaragoza, 50009, Spain. (e-mail: \texttt{rbarrio@unizar.es}).

Marcos Rodr\'{\i}guez is with CoDy group, Zaragoza, Universidad de Zaragoza, Pedro Cerbuna, 12, Zaragoza, 50009, Spain.
(e-mail: \texttt{imark.rodriguez@gmail.com}).}
% <-this % stops a space
\thanks{}% <-this % stops a space
\thanks{}}

\maketitle

\begin{abstract}
In image processing, \emph{Total Variation (TV) regularization models} are commonly used to recover blurred images. One of the most efficient and popular methods to solve the convex TV problem is the \emph{Alternating Direction Method of Multipliers} (ADMM) algorithm, recently extended using the inertial proximal point method. Although all the classical studies focus on only a convex formulation, recent articles are paying increasing attention to the nonconvex methodology  due to its good numerical performance and properties. In this paper, we propose to extend the classical formulation with a novel \emph{nonconvex Alternating Direction Method of Multipliers with the Inertial technique} (IADMM). Under certain assumptions on the parameters, we prove the convergence  of the algorithm with the help of the Kurdyka-{\L}ojasiewicz property.
We also present numerical simulations on classical TV image reconstruction problems  to illustrate the efficiency  of the new algorithm  and its behavior compared with the well established ADMM method.
\end{abstract}


\IEEEpeerreviewmaketitle

\section{Introduction}
Denoising and deblurring have numerous applications in communications, control, machine learning, and many other fields of engineering and
science. Restoration of distorted images is, from the theoretical, as well as from the practical point of view, one of the most interesting and important problems of image processing. One special case is the blurring, due, for instance, to incorrect focus and/or blurring due to movement, or added Gaussian noise (a Gaussian blur).

A mathematical model for the process of  blurred images can be expressed as follows. Let $\Omega \equiv \mathbb{R}^{N^2}$ be a two-dimensional index set representing the image domain, $\tilde{u}\in\Omega$
be the original image, $\tilde{f}$ be the observed image, and $\tilde{K}$ be a linear blurring operator. Then, the blurred image can be written \cite{Hansen:2006} as
\begin{eqnarray}
\tilde{f}=\tilde{K}(\tilde{u})+e,
\end{eqnarray}
where $e$ is an unknown additive noise vector. In this paper, the  blurring operator $\tilde{K}$ is assumed to be known, otherwise, one will deal with   the blind image deblurring problem \cite{lucy1974iterative}, in which  $\tilde{K}$   also needs to be solved.

In image processing one typically aims at recovering an image from noisy data while still
keeping edges in the image, and this goal is the main reason of the tremendous success of the \emph{Total Variation (TV) regularization} \cite{rudin1992nonlinear} for solving the deblurring problem (although other methods are also used). The TV method can be presented as
\begin{equation}\label{debblurring}
   \min_{\tilde{u}}\left(\frac{1}{2}\|\tilde{K}\tilde{u}-\tilde{f}\|^2_{L^2(\Omega)}+\sigma \|\nabla \tilde{u}\|_{L^1(\Omega^*)}\right),
\end{equation}
being $\sigma>0$ a parameter and where $\nabla$  is the gradient operator, $\Omega^*:=\nabla(\Omega)$, and the norms $\|\tilde{v}\|_{L^1(\Omega^*)}:=\sum_{i,j\in\Omega^*}|\tilde{v}_{i,j}|$, $\|\tilde{u}\|_{L^2(\Omega)}:=\sqrt{\sum_{i,j\in\Omega}|\tilde{u}_{i,j}|^2}$.

In most situations, rather than directly minimizing the
support of the image, one is interested in minimizing the support of the gradient of
the recovered image. In most references, the convex methodology is considered \cite{chen2015inertial,wang2008new,chen2018segmenting}, but in recent years, some nonconvex methods have been developed \cite{hintermuler2013nonconvex,hintermuller2014smoothing,xu2011image}. The use of a suitable
 nonconvex and nondifferentiable
function allows
possibly a smaller number of measurements than the convex one in compressed sensing \cite{hintermuler2013nonconvex}. In \cite{Nikolova2008} the authors showed that nonconvex regularization terms in
total variation-based image restoration yields even better edge preservation when compared
to the convex-type regularization. Moreover, they showed that it seems to be
also more robust with respect to noise.
Nonconvex regularization in image restoration poses significant challenges on the existence of solutions
of associated minimization problems and on the development of efficient solution algorithms.

The main  difference  between the convex and nonconvex methods is replacing the $l_1$-norm of the variational term by the nonconvex and nondifferentiable
function  $\|v\|_{\varphi(\Omega^*)}:=\sum_{i,j\in\Omega^*}\varphi(v_{i,j})$ that uses the nonconvex regulation function $\varphi$, and that we refer to as the $\ell^q$ semi-norm ($0<q<1$). Therefore, the general nonconvex deblurring model is presented as
\begin{equation}\label{debblurring2}
   \min_{\tilde{u}}\left(\frac{1}{2}\|\tilde{K}\tilde{u}-\tilde{f}\|^2_{L^2(\Omega)}+\sigma \|\nabla \tilde{u}\|_{\varphi(\Omega^*)}\right).
\end{equation}

Many efficient numerical
algorithms have been developed for solving the TV regularization problem.
One of the most efficient methods for the convex problem (\ref{debblurring}) is the Alternating Direction Method of Multipliers (ADMM) algorithm  \cite{boyd2011distributed,eckstein1994some,gabay1976dual}. In the general case (convex and nonconvex cases depend of the function $\varphi$), the method is constructed by introducing an auxiliary variable $\tilde{v}$, which actually represents $\nabla \tilde{u}$,  to reformulate \eqref{debblurring2}  into a composite optimization problem with linear constraints.
The augmented
Lagrange dual  function is   then
\begin{align}\label{lari-}
  &\mathcal{\tilde{L}}_{\delta}^{\varphi}(\tilde{u},\tilde{v},\tilde{p}):=\frac{1}{2}\|\tilde{K}\tilde{u}-\tilde{f}\|^2_{L^2(\Omega)}+\sigma \|\tilde{v}\|_{\varphi(\Omega^*)}\nonumber\\
   &\qquad\qquad-\langle \tilde{p},\nabla \tilde{u}-\tilde{v}\rangle+\frac{\delta}{2}\|\nabla  \tilde{u}-\tilde{v}\|_{L^2(\Omega^*)}^2,
\end{align}
where $\delta>0$ is a parameter, and the norm $\|X\|_{L^2(\Omega^*)}:=\sqrt{\sum_{i,j\in\Omega^*}|X_{i,j}|^2}$. If $\varphi(\cdot)=|\cdot|_1$, we use the notation $\mathcal{\tilde{L}}_{\delta}^{1}(\tilde{u},\tilde{v},\tilde{p})$ for representing $\mathcal{\tilde{L}}_{\delta}^{\varphi}(\tilde{u},\tilde{v},\tilde{p})$. Now, the standard convex ADMM method ($\varphi(\cdot)=|\cdot|_1$) for the deblurring problem can be presented as
 \begin{eqnarray}\label{oscheme}
 \left\{
   \begin{array}{lcl}
     \tilde{v}^{k+1}&=&\textrm{arg}\min_{\tilde{v}}\mathcal{\tilde{L}}_{\delta}^1(\tilde{u}^k,\tilde{v},\tilde{p}^k) \\[5.pt]
     \tilde{u}^{k+1}&=&\textrm{arg}\min_{\tilde{u}} \mathcal{\tilde{L}}_{\delta}^1(\tilde{u},\tilde{v}^{k+1},\tilde{p}^{k}) \\[5.pt]
     \tilde{p}^{k+1}&=&\tilde{p}^k-\delta(\nabla \tilde{u}^{k+1}-\tilde{v}^{k+1}) \\
   \end{array}
 \right.
\end{eqnarray}

The earlier analyses of convergence and performance of the ADMM algorithms directly depended on the existing results of ADMM framework \cite{boley2013local,deng2012global,he2012non,he20121,hong2012linear}. More recently, motivated by acceleration techniques proposed in \cite{polyak1964some}, inertial algorithms have been proposed for many areas such as (distributed) optimization and imaging sciences in references  \cite{banert2015backward,bot2015penalty,chen2015general,ochs2014ipiano,sun2018noner}. The ideas of the inertial strategy have been also applied to ADMM in \cite{chen2015inertial,bot2014inertial}; and under several  assumptions in convex case, some convergence results are proved in those articles. As the nonconvex penalty functions perform more efficiently in some applications, as above commented, nonconvex ADMM has been also developed and studied \cite{li2014splitting,hong2016convergence,sun2017iteratively,sun2017alternating,wang2014convergence,wang2015global,sun2017convergence,xu2016empirical,taylor2016training}.  The main goal of this paper is to propose a new algorithm that combines the nonconvex methods and the inertial strategy organically.

 In this paper, when $\varphi$ is nonconvex, we consider a new  inertial scheme for the image deblurring model (\ref{debblurring2}). One of the main differences (and new difficulties) with the convex ADMM, is that in order to properly define the nonconvex ADMM some extra assumptions are needed to prove the convergence. First, at least one of the objective functions has to be smooth. And more,  matrix corresponding to the smooth function is required to be injective, i.e., reversible. Thus, a direct application of the ADMM scheme to the image deblurring model cannot guarantee the convergence because the operator $\nabla$  fails to be injective (although the numerical performance may be  good in some cases). Considering this, we first modify the model (\ref{debblurring2}), and then we develop the new nonconvex inertial ADMM.
By using the  Kurdyka-{\L}ojasiewicz property, we prove the convergence of the new algorithm under several requirements on the parameters. In opposite to the convex case, selecting a suitable parameter $\delta$ is crucial to obtain the convergence of the new algorithm. In order to make the method more useful, we
provide a probabilistic strategy for selecting a suitable $\delta$.

The rest of the paper is organized as follows. In Section~\ref{sec:2} we collect some mathematical preliminaries needed for the convergence analysis. Section~\ref{sec:3} presents the details for the new algorithm (inertial alternating minimization algorithm, IADMM) including the schemes and parameters. In Section~\ref{sec:4}, we prove the convergence of the new algorithm.  Section~\ref{sec:5} reports the numerical results and compares the algorithm with convex and nonconvex ADMM. Section~\ref{sec:6} gives some conclusions. Finally, we provide in the Appendixes all the detailed proofs of the proposed results.

\section{Mathematical tools}
\label{sec:2}
In this section we present the definitions and basic properties of the subdifferentials
and the Kurdyka-{\L}ojasiewicz functions used later in the convergence analysis.
The basic notations  used in this paper are detailed in Table~\ref{tablebasic}.
\begin{table}
\begin{tabular}{c c}
  \hline
  \multicolumn{2}{c}{$\|\cdot\|$ stands for $\|\cdot\|_{L^2(\Omega)}$ or $\|\cdot\|_{L^2(\Omega^*)}$ ($L_2$ norm)}\\
    \hline
  \textrm{dist}$(x,C):=\min_{y\in C}\|x-y\|$  & $\|A\|_2:=\max_{\|x\|=1}\|Ax\|$    \\
  \hline
$\|\cdot\|_{\varphi}:=\|\cdot\|_{\varphi(\Omega^*)}$ & $\otimes$ stands for the Kronecker product\\
  \hline
   \multicolumn{2}{c}{$\mathcal{C}^1$ stands for the function class whose derivatives are continuous} \\
  \hline
\multicolumn{2}{c}{for a matrix A, \textrm{rank}(A) stands for its rank, $\textrm{Null}(A):=\{x\mid Ax=\textbf{0}\}$}\\
\hline
\multicolumn{2}{c}{$\lambda_{\min}(A)$ stands for the minimum eigenvalue of $A$}\\
\hline
\end{tabular}
\caption{Basic notations, where $x,y$ stand for points, $C$ stands for a set, $A$ stands for matrix \label{tablebasic}}
\end{table}
\subsection{Subdifferentials}
We collect several definitions as well as some useful properties in variational and convex
analysis (see the monographes \cite{mordukhovich2006variational,nesterov2004introductory,rockafellar2009variational,rockafellar2015convex}).
 For any matrix $A$, we define $A^*$ to be the adjoint of $A$.

\begin{definition} Let  $J: \mathbb{R}^N \rightarrow (-\infty, +\infty]$ be a proper and lower semicontinuous function. The (limiting) subdifferential, or simply the subdifferential, of $J$ at $x\in \mathbb{R}^N$, written as $\partial J(x)$, is defined  as
\begin{align*}
    \partial J(x):=\{u\in\mathbb{R}^N: \exists~  x^k\rightarrow x,~u^k\rightarrow u,~\textrm{such~ that}\\
    \lim_{y\neq x^k}\inf_{y\rightarrow x^k}\frac{J(y)-J(x^k)-\langle u^k, y-x^k\rangle}{\|y-x^k\|_2}\geq 0~\}.
\end{align*}
\end{definition}

It is easy to verify that the Fr\'echet subdifferential is convex and closed while the subdifferential is closed. When $J$ is convex,  the definition agrees with the subgradient in convex analysis \cite{rockafellar2015convex}.
We call $J$ is \emph{strongly convex} with constant $a \in \mathbb{R}^+$ if for any $x,y\in \textrm{dom}(J)$ and any $v\in \partial J(x)$, it holds that
$J(y)\geq J(x)+\langle v,y-x\rangle+\frac{a}{2}\|y-x\|_2^2.$
And $J$ is called as  $L$-\emph{gradient continuous} (Lipschitz) with constant $L>0$ if
$\|\nabla J(x)-\nabla J(y)\|_2\leq L\|x-y\|_2.$
Noting the closedness of the subdifferential, we have the following simple proposition.

\smallskip

\begin{proposition}\label{sublimit}
If $v^k\in \partial J(x^k)$, $\lim_{k}v^k=v$ and $\lim_{k}x^k=x$, then we have
$
    v\in \partial J(x).
$
\end{proposition}

\begin{definition}
A necessary condition for $x\in\mathbb{R}^N$ to be a minimizer of $J(x)$ is
\begin{equation}\label{Fermat}
\textbf{0}\in \partial J(x),
\end{equation}
which is also sufficient
when $J$ is convex.
A point that satisfies (\ref{Fermat}) is called \emph{(limiting) critical point}. The set of critical points of $J(x)$ is denoted by $\textrm{crit}(J)$.
\end{definition}
\smallskip

With these basics, we can easily obtain the following proposition.

\begin{proposition}\label{criL}
If  $(\bar{u},\bar{v},\bar{p})$ is a critical point of $\mathcal{L}^{\varphi}_{\delta}$ whose definition given in \eqref{lari}, it must hold that
\begin{eqnarray}
\left\{\begin{array}{lcl}
         -\bar{p}&\in&\sigma\partial \|\bar{v}\|_{\varphi}, \\[5.pt]
         \nabla^{*}\bar{p}&=& K^{*}(K\bar{u}-f), \\[5.pt]
         \nabla \bar{u}&=& \bar{v}.
       \end{array}
\right.
\end{eqnarray}
\end{proposition}
\smallskip

Finally, the \emph{proximal map} of $J$ is defined as
\begin{eqnarray}
\label{proximal}
\mathcal{S}_{J}(x)\in\textrm{arg}\min_{y}\bigg\{J(y)+\frac{1}{2}\|x-y\|_2^2\bigg\}.
\end{eqnarray}
Note that $J$ can be nonconvex. If $J$ is convex, $\mathcal{S}_{J}$ is a point-to-point operator; otherwise, it may be  point-to-set.

\subsection{Kurdyka-{\L}ojasiewicz property}

In this paper the convergence analysis is based on the Kurdyka-{\L}ojasiewicz functions, originated in the seminal works of {\L}ojasiewicz \cite{Lojasiewicz63} and Kurdyka \cite{Kurdyka98}.  This kind of functions has played a key role in several recent convergence results on nonconvex minimization problems and they are ubiquitous in applications.

\begin{definition}[\cite{attouch2013convergence,bolte2014proximal}]\label{KL}
(a) The function $J: \mathbb{R}^N \rightarrow (-\infty, +\infty]$ is said to have the Kurdyka-{\L}ojasiewicz property at $\widehat{x}\in \textrm{dom}(\partial J)$ if there
 exist $\eta\in (0, +\infty]$, a neighborhood $U$ of $\widehat{x}$ and a continuous concave function $\rho: [0, \eta)\rightarrow \mathbb{R}^+$ such that
\begin{enumerate}
  \item $\rho(0)=0$.
  \item $\rho$ is $\mathcal{C}^1$ on $(0, \eta)$.
  \item For all $s\in(0, \eta)$, $\rho'(s)>0$.
  \item For all $x$ in $U\bigcap\{x|J(\widehat{x})<J(x)<J(\widehat{x})+\eta\}$, the Kurdyka-{\L}ojasiewicz inequality holds
\begin{equation}
  \rho'(J(x)-J(\widehat{x}))\cdot\textrm{dist}({\rm{\bf{0}}},\partial J(x))\geq 1.
\end{equation}
\end{enumerate}

(b) Proper lower semicontinuous functions which satisfy the Kurdyka-{\L}ojasiewicz inequality at each point of $\textrm{dom}(\partial J)$ are called \emph{K{\L} functions}.
\end{definition}
\smallskip

\begin{remark} There are large sets of functions that are K{\L} functions \cite{attouch2013convergence}.
\end{remark}

\begin{lemma}[\cite{bolte2014proximal}]\label{con}
Let $J: \mathbb{R}^N\rightarrow \mathbb{R}$ be a proper lower semi-continuous function and $\Pi$ be a compact set. If $J$ is a constant on $\Pi$ and $J$ satisfies the K{\L} property at each point on $\Pi$, then there exists a concave function $\rho$ satisfying the four properties given in Definition \ref{KL}, and constants $\eta,\varepsilon>0$ such that for any $\widehat{x}\in \Pi$ and any $x$ satisfying that $\textrm{dist}(x,\Pi)<\varepsilon$ and $f(\widehat{x})<f(x)<f(\widehat{x})+\eta$, it holds that
\begin{equation}
    \rho'(J(x)-J(\widehat{x}))\cdot\textrm{dist}({\rm{\bf{0}}},\partial J(x))\geq 1.
\end{equation}
\end{lemma}

\section{Nonconvex IADMM algorithm}
\label{sec:3}
In this section we introduce the new extended Inertial Alternating Direction Method of Multipliers (ADMM) algorithm for nonconvex functions.

In this  paper, we consider $\Omega=[1,2,\ldots,N]\times [1,2,\ldots,N]$ (equivalent to space $\mathbb{R}^{N^2}$) as the two-dimensional index set representing the image domain.  In this case, the image variable constrained on $\Omega$ is actually a $N\times N$ matrix.  We use the symbol $\tilde{u}$ to present its vectorization (a vector of all the columns of the image variable).
And then the original total variation operator then enjoys the following  form
\begin{eqnarray}
\nabla (\tilde{u})=\left(
                 \begin{array}{c}
                   I_N\otimes D\\
                   D\otimes I_N \\
                 \end{array}
               \right)\cdot \tilde{u},
\end{eqnarray}
where  $I_N$ the identity matrix and $D$ the banded matrix $$D=\left(
                                                          \begin{array}{ccccc}
                                                            1 & -1 &  &  &\\
                                                             & 1 & \ddots & & \\
                                                             &  & \ddots & -1& \\
                                                             &  &  & 1&   -1 \\
                                                          \end{array}
                                                        \right)_{(N-1)\times N}.
$$
If we directly apply  the inertial ADMM, the convergence is hard to be proved as $\nabla$ fails to be injective. Therefore, we need to modify the image deblurring model (\ref{debblurring2}). To that goal, we define
\[
\mathcal{T}:=\left(
                     \begin{array}{cc}
                     I_N\otimes D &  \\
                        &          D\otimes I_N \\
                     \end{array}
                   \right), \quad u:=\left(
                                                                \begin{array}{c}
                                                                  u_1\in \mathbb{R}^{N^2} \\
                                                                  u_2 \in \mathbb{R}^{N^2}\\
                                                                \end{array}
                                                              \right).
\]
Obviously, we have $\mathcal{T}\in \mathbb{R}^{2N(N-1)\times 2N^2}$.
Noting that
\[
\begin{array}{l}
\textrm{rank}(\mathcal{T})=\textrm{rank}(I_N\otimes D)+\textrm{rank}(D\otimes I_N)\\
=\textrm{rank}(I_N)\cdot\textrm{rank}(D)+\textrm{rank}(D)\cdot\textrm{rank}(I_N)=2N(N-1),
\end{array}
\]
and thus, $\mathcal{T}$ is injective.
The following technical lemma focuses on giving a lower bound for the operator $\mathcal{T}^*$.

\begin{lemma}
\label{lemmaini}
For any $x\in \mathbb{R}^{2N(N-1)}$, it holds that
\begin{eqnarray}
\|\mathcal{T}^{*}x\|\geq\frac{\|x\|}{\theta},
\end{eqnarray}
where $\theta=1/(2\sin(\frac{\pi}{2N})\big)$.
\end{lemma}
\smallskip

Then, the image deblurring model (\ref{debblurring2}) is equivalent to
 \begin{equation}\label{debblurring3}
      \min_{u_1=u_2}\left(\frac{1}{2}\|\tilde{K}u_1-\tilde{f}\|^2+\sigma \|\mathcal{T} u\|_{\varphi}\right).
 \end{equation}
 Instead, we consider its extended penalty form
  \begin{equation}\label{debblurring4}
      \min_{u}\left(\frac{1}{2}\|\tilde{K}u_1-\tilde{f}\|^2+\sigma \|\mathcal{T} u\|_{\varphi}+\frac{\beta^2}{2}\|u_1-u_2\|^2\right),
 \end{equation}
 where $\beta>0$ is a large weight parameter. Therefore, we apply the nonconvex  inertial ADMM to
   \begin{equation}\label{debblurring-final}
      \min_{u\in \mathbb{R}^{2N^2}}\left(\frac{1}{2}\|Ku-f\|^2+\sigma \|\mathcal{T} u\|_{\varphi}\right),
 \end{equation}
 where $K=\left(
            \begin{array}{cc}
              \tilde{K} & \textbf{0} \\
              \beta \mathbb{I} & -\beta \mathbb{I} \\
            \end{array}
          \right)
 $ and $f=\left(
            \begin{array}{c}
              \tilde{f} \\
              \textbf{0} \\
            \end{array}
          \right)
 $.
 This leads us to define the function
\begin{align}\label{lari}
   &\mathcal{L}_{\delta}^{\varphi}(u,v,p):=\frac{1}{2}\|Ku-f\|^2\nonumber\\
   &\quad\quad+\sigma \|v\|_{\varphi}-\langle p,\mathcal{T} u-v\rangle+\frac{\delta}{2}\|\mathcal{T} u-v\|^2.
\end{align}

Inertial methods have witnessed  great success in convex ADMM and nonconvex first-order algorithms. In the nonconvex optimization community, the inertial style ADMM has never been proposed and analyzed.  The convex inertial  ADMM has been proposed in \cite{chen2015inertial}, in which one  first uses the ``inertial method" to refresh the current sequence with last iteration, and then performs the ADMM scheme with the updated variables. However, the direct extension of convex ADMM   is not allowed   in the nonconvex settings.  This is because without convexity, several descents are heavily dependent on the continuities of the functions, which $\varphi$ may fail to obey.  And the difference of function values at two  different   points is hard to estimate, which leads to  troubles in the convergence proof. Thus, in the updating of $v^{k+1}$, we used $u^k$ rather than the updated one. And the nonconvex IADMM scheme proposed in this paper is defined as follows
\begin{eqnarray}
\label{debeq}
\left\{\begin{array}{lcl}
  (\widehat{u}^k, \widehat{v}^k, \widehat{p}^k) &=& (u^k, v^k, p^k)\\[4.pt]
  && \hspace*{-1.cm}+\alpha(u^k-u^{k-1}, v^k-v^{k-1}, p^k-p^{k-1}), \\[4.pt]
  v^{k+1}&=&\textrm{arg}\min_{v} \mathcal{L}_{\delta}^{\varphi}(u^k,v,\widehat{p}^k), \\[4.pt]
  u^{k+1}&=&\textrm{arg}\min_{u} \mathcal{L}_{\delta}^{\varphi}(u,v^{k+1},\widehat{p}^k), \\[4.pt]
  p^{k+1}&=&\widehat{p}^k-\delta(\mathcal{T} u^{k+1}-v^{k+1}),
\end{array}\right.
\end{eqnarray}
where $\alpha>0$ is a free parameter chosen by the user.   Actually, if $\alpha=0$, the algorithm then will reduce to basic ADMM.

Now we can focus on rewriting the inertial scheme for the image deblurring model (\ref{debeq}).
First, we rearrange the minimization for $v^{k+1}$,
\begin{eqnarray}
\label{vk1}
v^{k+1}=\mathcal{S}_{\frac{\sigma}{\delta}\|\cdot\|_{\varphi(\Omega^*)}}\bigg(\mathcal{T} u^k-\frac{\widehat{p}^k}{\delta}\bigg),
\end{eqnarray}
being $S_J$ the proximal map of $J$ (\ref{proximal}).
For a matrix $v$, and indices $(i,j)\in \Omega^*$,
\begin{eqnarray}
[\mathcal{S}_{\frac{\sigma}{\delta}\|\cdot\|_{\varphi}}(v)]_{i,j}=
\mathcal{S}_{\frac{\sigma}{\delta}\varphi}(v_{i,j}).
\end{eqnarray}
The scheme for updating $u^{k+1}$ can be rewritten as
\begin{equation}
    K^*(Ku^{k+1}-f)-\mathcal{T}^* \widehat{p}^k+\delta\mathcal{T}^*(\mathcal{T} u^{k+1}-v^{k+1})=\textbf{0}.
\end{equation}
That is also
\begin{equation}
\label{uk1}
    u^{k+1}=(K^*K+\delta\mathcal{T}^*\mathcal{T})^{-1}(K^* f+\delta\mathcal{T}^*v^{k+1}+\mathcal{T}^* \widehat{p}^k).
\end{equation}

Taking into account (\ref{debeq}), (\ref{vk1}) and (\ref{uk1}) we propose a nonconvex inertial version of the IADMM algorithm (Algorithm~\ref{algo}).

\begin{algorithm}[h]
\caption{Nonconvex Inertial Alternating Minimization (IADMM) for Image Deblurring \label{algo}}
\begin{algorithmic}
\REQUIRE parameters $\alpha>0, \, \delta>0,$\\
\textbf{Initialization}: $u^0=u^1,\,v^{0}=v^1,\,p^0=p^1$\\
\textbf{for}~$k=1,2,\ldots$ \\
~~~ $(\widehat{u}^k, \widehat{v}^k, \widehat{p}^k) = (u^k, \, v^k, \, p^k)$\\[4.pt]
~~~~~~~~~~ $+\alpha(u^k-u^{k-1}, \, v^k-v^{k-1}, p^k-p^{k-1})$\\[7.pt]
~~~ $v^{k+1}=\mathcal{S}_{\frac{\sigma}{\delta}\|\cdot\|_{\varphi}}\bigg(\displaystyle\mathcal{T} u^k-\frac{\widehat{p}^k}{\delta}\bigg)$ \\[7.pt]
~~~ $u^{k+1}=(K^*K+\delta\mathcal{T}^*\mathcal{T})^{-1}(K^* f+\delta\mathcal{T}^*v^{k+1}+\mathcal{T}^*\widehat{p}^k)$ \\[5.pt]
~~~ $ p^{k+1}=\widehat{p}^k-\delta(\mathcal{T} u^{k+1}-v^{k+1})$ \\
\textbf{end for}\\
\end{algorithmic}
\end{algorithm}

\textbf{Assumption 1:}  We assume that $\textrm{Null}(K_{0})\bigcap \textrm{Null}(\mathcal{T})=\{\textbf{0}\}$, where $K_0:=\left(
            \begin{array}{cc}
              \tilde{K} & \textbf{0} \\
               \mathbb{I} & - \mathbb{I} \\
            \end{array}
          \right)$. And the minimum single value is given as $\nu$. \\
          \smallskip

This hypothesis also indicates that the matrix  $\left(
     \begin{array}{c}
       \mathcal{T} \\
       K_0\\
     \end{array}
   \right)$ is reversible.
Note that the rank of $\mathcal{T}$ is $2N(N-1)$.   Then, the assumed hypothesis is easy to be satisfied.

We remark that $u^{k+1}$ is the minimizer of $ \mathcal{L}_{\delta}^{\varphi}(u,v^{k+1},\widehat{p}^k)$, which is strongly convex with $\lambda_{\min}(K^*K+\delta\mathcal{T}^*\mathcal{T})$. If we set $\delta,\rho>1$, then we have
\small
\begin{align}\label{des+}
 &\mathcal{L}_{\delta}^{\varphi}(u^{k},v^{k+1},\widehat{p}^k)- \mathcal{L}_{\delta}^{\varphi}(u^{k+1},v^{k+1},\widehat{p}^k)\nonumber\\
 &\geq \frac{\lambda_{\min}(K^*K+\delta\mathcal{T}^*\mathcal{T})}{2}\|u^{k+1}-u^{k}\|^2\geq\frac{\nu}{2}\|u^{k+1}-u^{k}\|^2,
\end{align}
\normalsize
where we used the fact $\lambda_{\min}(K^*K+\delta\mathcal{T}^*\mathcal{T})\geq \lambda_{\min}(K^*_0K_0+\mathcal{T}^*\mathcal{T})$ when $\delta,\beta\geq 1$.

The following problem is what $\nu$ exactly is. In a real situation, the dimensions of $\mathcal{T}$ are large, and so, a direct calculation leads to a large computational cost. Therefore, we provide a probabilistic method to estimate a suitable value of $\theta$. If $K^*K+\mathcal{T}^*\mathcal{T}$ is reversible, it is easy to see that
${1}/{\nu}=\|(K^*K+\mathcal{T}^*\mathcal{T})^{-1}\|_2$. Then, if we obtain a bound  $\|(K^*K+\mathcal{T}^*\mathcal{T})^{-1}\|_2\leq c$, we then have $\nu\geq{1}/{c}$. To that goal we employ a Lemma proposed in \cite{halko2011finding}:
\begin{lemma}[Lemma 4.1, \cite{halko2011finding}]
For a fixed positive integer $M$ and a
real number $b>1$, and given an independent family $\{w_i\}_{i=1,2,\ldots,M}$ of standard
Gaussian vectors, we have that
$$\|(K^*K+\mathcal{T}^*\mathcal{T})^{-1}\|_2\nonumber\\
\leq b\sqrt{\frac{2}{\pi}}\max_{i=1,2,\ldots,M}\|(K^*K+\mathcal{T}^*\mathcal{T})^{-1}w_i\|_2
$$
with probability at least $1-b^{-M}$.
\end{lemma}
\smallskip

Note that for computing $(K^*K+\mathcal{T}^*\mathcal{T})^{-1}w_i$ we just need several FFT and inverse FFT. Therefore, its  computational cost is low ($\mathcal{O}(N \log N$), and the estimation of $\nu$ is very fast.

\section{Convergence analysis}
\label{sec:4}
This section consists of two parts and provides a complete analysis of the convergence problem of the nonconvex IADMM algorithm.    The first subsection contains the main convergence results, the proof sketch, the difficulties in the proof and theoretical  contributions. While the second subsection introduces the necessary technical lemmas.
Assumption 1 holds through this section.

\subsection{Main results}
\begin{theorem}[Stationary point convergence]\label{th1}
Assume that the free parameter $\delta$ satisfies the condition
\begin{equation}\label{delta1}
 \delta>\max\left\{1, \, \frac{6\theta^2\|K\|_2^4}{\nu}+\frac{7\alpha^2\theta^2\|K\|_2^4}{\nu}\right\}
\end{equation}
with $\theta=1/(2\sin(\frac{\pi}{2N})\big)$, then
any cluster point $(u^*,v^*,p^*)$ is also a critical point of $\mathcal{L}^{\varphi}_{\delta}$.
\end{theorem}
\smallskip

Theorem \ref{th1} describes the stationary point convergence result for the IADMM method, which is free of using the K{\L} property of the functions.  If the K{\L} property is further assumed, the sequence convergence can be proved giving us the Theorem \ref{th2}.

\begin{theorem}[Sequence convergence]\label{th2}
 Let   condition~\eqref{condelta} hold, and the auxiliary function $F$ (\ref{fun}) be a K{\L} function. Then the sequence $\{(u^k,v^k,p^k)\}_{k=0,1,2,\ldots}$, generated by Algorithm~\ref{algo}, converges to a critical point of $\mathcal{L}^{\varphi}_{\delta}$.
\end{theorem}
\smallskip

The proof  can be divided into two parts, and in order to help the reader we first give a \emph{brief sketch of the proof}:

\textbf{I}. In the first part we introduce  an auxiliary sequence $\{w^{k}\}_{k=0,1,2,\ldots}$, where  $w^k$  are composite points from $\{(u^k,v^k,p^k)\}_{k=0,1,2,\ldots}$.  An auxiliary function $F$ is also proposed. In Lemma~\ref{descend1}, we prove a ``sufficient descent condition'' of the values of $F$ at
$\{w^{k}\}_{k=0,1,2,\ldots}$, i.e.,
\begin{equation}
    F(w^k)-F(w^{k+1})\geq \widehat{h}\|u^{k+1}-u^k\|^2,
\end{equation}
where $\widehat{h}>0$ is a positive constant.

\textbf{II}. We prove a ``relative error condition" of
$\{w^{k}\}_{k=0,1,2,\ldots}$, i.e., there exists $v^{k+1}\in \partial F(w^{k+1})$ such that
\begin{equation}
    \|v^{k+1}\|\leq\gamma(\|u^{k+1}-u^k\|+\|u^{k}-u^{k-1}\|),
\end{equation}
where $\gamma>0$ is a positive constant.  Note that this condition is different from the ``real'' relative error condition proposed in \cite{attouch2013convergence}.

The major  difficulty in  deriving these two conditions is the use of inertial terms, with which  the descent values are lower bounded  by  $\|v^{k+1}-\widehat{v}^k\|^2$ and $\|p^{k+1}-\widehat{p}^k\|^2$ rather than  $\|v^{k+1}-v^k\|^2$ and $\|p^{k+1}-p^k\|^2$. Similarly, the  relative error is bounded by $\|v^{k+1}-\widehat{v}^k\|$ and $\|p^{k+1}-\widehat{p}^k\|$. The relative error can be expanded by triangle inequalities, which is relatively proved. However, for the sufficient descent, the use of the triangle inequalities is much more difficult and technical for the lower boundedness.

The theoretical contributions in this paper are two-fold. The first one is, of course, dealing with the difference caused by the inertial terms. This part also includes how to design the scheme of the algorithm, whose details have been presented in previous section.   The second one is to determine the parameters for IADMM  applied to the image deblurring.

\smallskip
The main results  can be proved with the following  lemmas.
\smallskip

\subsection{Technical lemmas}

\begin{lemma}\label{fp}
Let the sequence $\{(u^k,v^k,p^k)\}_{k=0,1,2,\ldots}$ be generated by Algorithm~\ref{algo} to solve problem~(\ref{debblurring2}), then
\begin{equation}\label{condition1}
     \|p^k-p^{k+1}\|\leq \theta\|K\|_2^2 \, \|u^{k+1}-u^{k}\|,
\end{equation}
where $\theta$ is given in Lemma \ref{lemmaini} and $\|K\|_2$ denotes the spectral radius of $K$.
\end{lemma}
\smallskip

Now we provide the main technical lemma that states the descent condition for a suitable function of the sequences of the Algorithm~\ref{algo}.

\begin{lemma}\label{descend1}
Let the sequence $\{(u^k,v^k,p^k)\}_{k=0,1,2,\ldots}$ be generated by Algorithm~\ref{algo} and conditions of Theorem \ref{th1} hold. By defining the auxiliary function $F$
\begin{equation}\label{fun}
    F(u,v,p,x):=\mathcal{L}_{\delta}^{\varphi}(u,v,p)+
    \displaystyle\frac{7\alpha^2\theta^2\|K\|_2^4}{2\delta}\|u-x\|^2,
\end{equation}
and
\begin{equation}
    w^k:=(u^{k},v^{k},p^{k},u^{k-1}),
\end{equation}
we have that
\begin{equation}
\label{descentcond}
    F(w^k)-F(w^{k+1})\geq \widehat{h}\|u^{k+1}-u^k\|^2,
\end{equation}
where
\begin{equation}
\label{auxh}\widehat{h}:=\frac{\nu}{2}-\bigg( \frac{3\theta^2\|K\|_2^4}{\delta}+\frac{7\alpha^2\theta^2\|K\|_2^4}{2\delta}\bigg)>0.\end{equation}
If the sequence  $\{(u^k,v^k,p^k)\}_{k=0,1,2,\ldots}$ is bounded, then, it holds that
\small
\begin{equation}
\lim_{k}\|u^{k+1}-u^k\|=\lim_{k}\|v^{k+1}-v^k\|=\lim_{k}\|p^{k+1}-p^k\|=0.
\end{equation}
\normalsize
\end{lemma}
\smallskip

\begin{remark}
Based on Lemma~\ref{descend1}, it is important to guarantee that the condition (\ref{delta1}) can be satisfied.
This fact can be reached if $\delta$ is large enough.  Fortunately, in the Algorithm~\ref{algo} the parameter $\delta$ can be fixed by the user. Thus, the parameter $\delta$ shall be chosen large enough to guarantee the convergence considering condition (\ref{delta1}).
\end{remark}
\smallskip

\begin{lemma}\label{bound}
If the nonconvex regulation function $\varphi$ is coercive and
\begin{equation}
\label{delta22}
\delta>\theta^2\|K\|_2^4,
\end{equation}
the sequence $\{w^k\}_{k=0,1,2,\ldots}$ is bounded.
\end{lemma}
\smallskip

\begin{remark}
By combining the  conditions (\ref{delta1}) and (\ref{delta22}), we just need
that
\begin{equation}\label{condelta}
    \delta>\max\left\{1,\,\frac{6\theta^2\|K\|_2^4}{\nu}+
    \frac{7\alpha^2\theta^2\|K\|_2^4}{\nu},\,\theta^2\|K\|_2^4\right\}.
\end{equation}
\end{remark}
\smallskip

\begin{lemma}\label{relative}
Let  the sequence $\{(u^k,v^k,p^k)\}_{k=0,1,2,\ldots}$ be generated by Algorithm~\ref{algo}. Then, for any $k$, there exists $\gamma>0$ and $s^{k+1}\in \partial F(w^{k+1})$ such that
\begin{equation}
    \|s^{k+1}\|\leq \gamma(\|u^{k+1}-u^{k}\|+\|u^{k}-u^{k-1}\|).
\end{equation}
\end{lemma}
\smallskip

Now, we recall a  definition about the  limit point set $\mathcal{M}$  introduced in \cite{bolte2014proximal}, which
denotes the set of all the stationary points generated by the nonconvex IADMM. The specific mathematical definition of $\mathcal{M}$ is given as follows.

\begin{definition}
\label{deflim}
Let $\{w^k\}_{k=0,1,2,\ldots}$ be generated by the nonconvex IADMM .
We define the set $\mathcal{M}$ by
\begin{align}
&\mathcal{M}:=\big\{x\in \mathbb{R}^N: \exists~\textrm{an increasing sequence of integers}\nonumber\\
& \qquad \qquad \{k_j\}_{j\in\mathbb{N}}~
\textrm{such that} ~x^{k_j}\rightarrow x ~\textrm{as}~ j\rightarrow \infty\big\}.
\end{align}
\end{definition}
\smallskip

\begin{lemma}\label{station}
Let the sequence $\{(u^k,v^k,p^k)\}_{k=0,1,2,\ldots}$ be generated by Algorithm~\ref{algo}, $F$ the auxiliary function defined in (\ref{fun}) and suppose that  condition (\ref{condelta}) holds. Then, we have the following results.

(1) $\mathcal{M}$ is nonempty and $\mathcal{M}\subseteq \textrm{crit}(F)$.

(2) $\lim_{k}\textrm{dist}(w^k,\mathcal{M})=0$.

(3) The objective function $F$ is finite and constant on $\mathcal{M}$.
\end{lemma}
\smallskip

\section{Numerics}
\label{sec:5}
In this section, we illustrate the effectiveness of the proposed algorithm on different numerical blurred images with Gaussian blur.

\begin{figure}[h]
\centering
 \includegraphics[width=0.9\columnwidth]{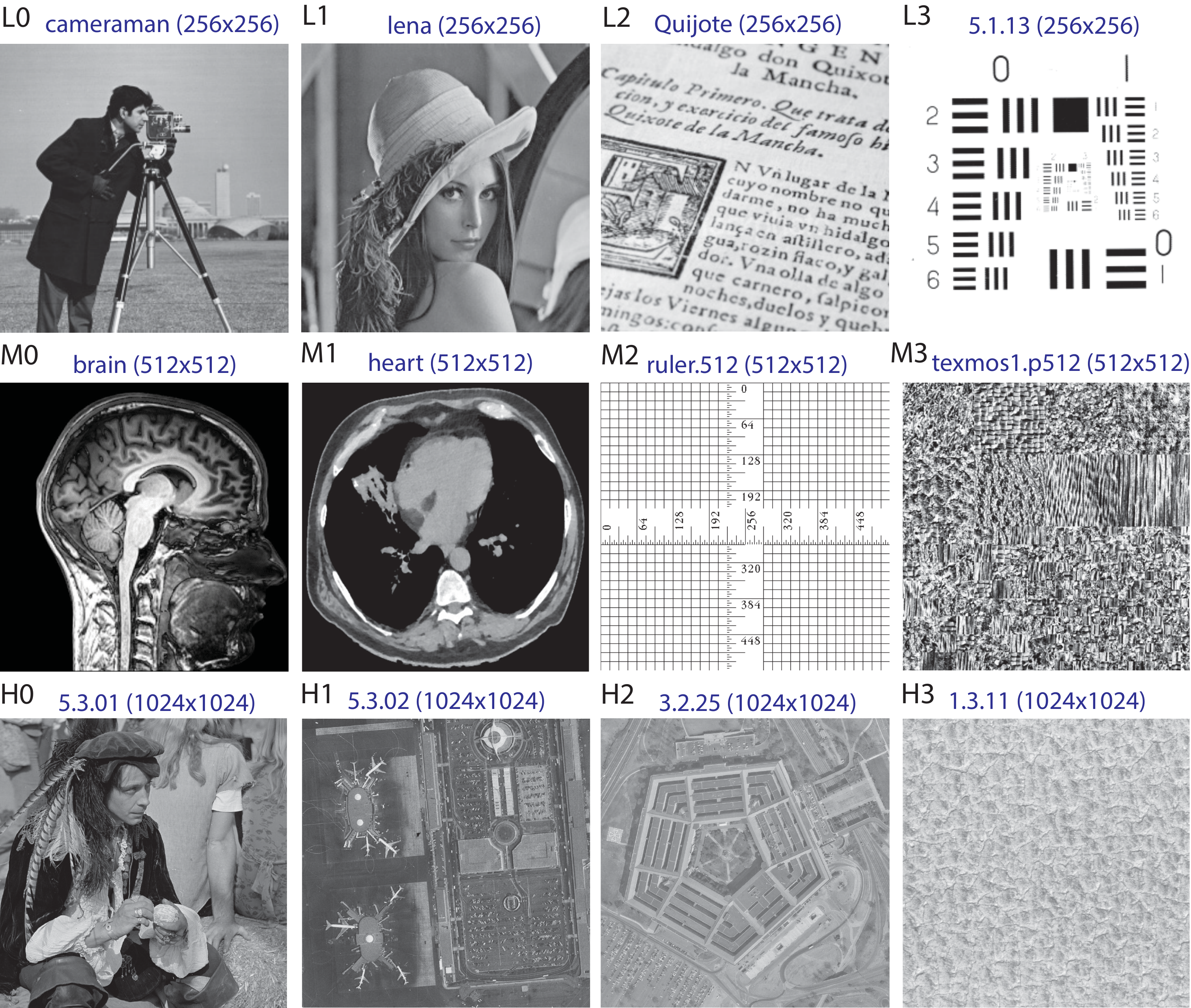}
  \caption{Original images (some of them from USC-SIPI image database) in low (L), medium (M) and hight-resolution (H). (L0) Cameraman  ($256\times256$); (L1) Lena ($256\times256$); (L2) ``El Quijote" ($256\times256$); (L3) 5.1.09 ($256\times256$); (M0) Brain ($512\times512$);
(M1) Heart ($512\times512$); (M2) ruler.512 ($512\times512$); (M3) texmos1.p512 ($512\times512$); (H0) 5.3.01 ($1024\times1024$);
(H1) 5.3.02 ($1024\times1024$); (H2) 3.2.25 ($1024\times1024$); (H3) 1.3.11 ($1024\times1024$). \label{fig1}}
\end{figure}

\begin{figure}[htbp]
\centering
 \includegraphics[width=0.85\columnwidth]{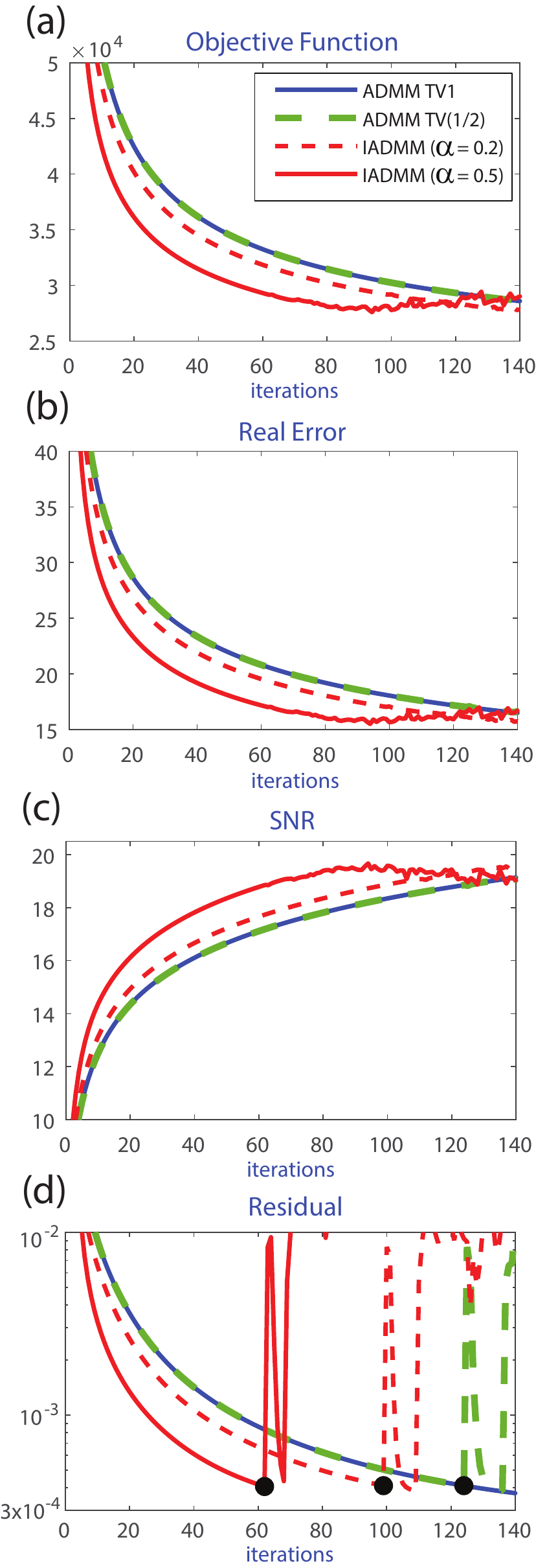}
  \caption{The evolution curves of the (a) \emph{objective function}, (b) the \emph{Real Error}, (c) the \emph{signal-to-noise ratio (SNR)} and (d) the \emph{Residual}, for the brain M0 image versus the iteration number. All the simulations have been done using the ADMM TV1 and TV(1/2) algorithms and the IADMM algorithm using $\alpha = 0.2$ and $0.5$, and $\delta = 0.001$.}
  \label{tests}
\end{figure}

All the programs have been written entirely in \texttt{C++}, and all the experiments are
implemented under Linux running on a desktop computer with an Intel Core
i5-2400S CPU (2.5 GHz) and 4 GB Memory. The FFT subroutines used in the algorithms are taken from the \texttt{fftw-3}\footnote{\texttt{http://www.fftw.org/}} library.
As test problems we have selected  twelve images (see Figure \ref{fig1}), which include seven images from USC-SIPI image database\footnote{\texttt{http://sipi.usc.edu/database/}}, two classical test images (Lena and cameraman), one text image from ``El Quijote" book and two medical images. In order to obtain the blurred images, we use,  as it is common in literature, the blurring operator generated using a convolution with Gaussian kernel (\texttt{KernelSize} $= 17 \times 17$, \texttt{KernelMu} $=0$, \texttt{KernelSigma} $= 7$) and circular mapping on the edges of the image.

\begin{table*}[htb!]
  \footnotesize
  \begin{center}{ deblurring results for $\delta = 0.001$, $\varepsilon = 0.005$}\end{center}
  \hspace*{-0.7cm}{\tt \noindent\begin{tabular}{|r|rrrr|rrrrr|rrrrr|}
    \cline{2-15}
    \multicolumn{1}{c|}{} & \multicolumn{4}{c|}{IADMM ($\alpha = 0.5$)}
        & \multicolumn{5}{c|}{IADMM ($\alpha = 0.2$)} & \multicolumn{5}{c|}{ADMM(TV1)}\\
    \hline
    \rowcolor[gray]{0.9}
    IMG & I1 & ERROR & SNR & RES & I2 & ERROR & SNR & RES & I2/I1 & I3 & ERROR & SNR & RES & I3/I1 \\
L0 &    4 &     17.4 & 11.1 & 4.59e-03 &    6 &     17.4 & 11.1 & 4.62e-03 & 1.50 &    7 &     17.6 & 11.0 & 4.80e-03 & 1.75 \\
\rowcolor[gray]{0.9}
L1 &    5 &     16.7 & 10.0 & 4.67e-03 &    8 &     16.5 & 10.1 & 4.60e-03 & 1.60 &    9 &     16.8 &  9.9 & 4.95e-03 & 1.80 \\
L2 &    2 &     19.4 &  8.7 & 4.36e-03 &    2 &     19.9 &  8.5 & 4.69e-03 & 1.00 &    2 &     20.2 &  8.3 & 4.88e-03 & 1.00 \\
\rowcolor[gray]{0.9}
L3 &    2 &     16.6 & 13.2 & 3.01e-03 &    2 &     17.1 & 13.0 & 3.25e-03 & 1.00 &    2 &     17.9 & 12.6 & 3.67e-03 & 1.00 \\
M0 &    8 &     30.6 & 13.8 & 4.50e-03 &   13 &     30.4 & 13.8 & 4.67e-03 & 1.62 &   16 &     30.6 & 13.8 & 4.88e-03 & 2.00 \\
\rowcolor[gray]{0.9}
M1 &    9 &     28.5 & 13.7 & 4.33e-03 &   14 &     28.7 & 13.6 & 4.86e-03 & 1.56 &   18 &     28.5 & 13.7 & 4.78e-03 & 2.00 \\
M2 &    6 &     36.6 & 12.8 & 3.48e-03 &    7 &     40.1 & 12.0 & 4.82e-03 & 1.17 &    9 &     39.9 & 12.1 & 4.86e-03 & 1.50 \\
\rowcolor[gray]{0.9}
M3 &    7 &     35.5 & 12.5 & 4.34e-03 &   11 &     35.6 & 12.4 & 4.50e-03 & 1.57 &   13 &     36.3 & 12.3 & 4.97e-03 & 1.86 \\
H0 &    8 &     71.4 & 10.3 & 4.50e-03 &   13 &     71.0 & 10.3 & 4.65e-03 & 1.62 &   16 &     71.4 & 10.3 & 4.84e-03 & 2.00 \\
\rowcolor[gray]{0.9}
H1 &    9 &     69.2 &  6.3 & 4.47e-03 &   14 &     69.7 &  6.2 & 4.96e-03 & 1.56 &   18 &     69.2 &  6.3 & 4.89e-03 & 2.00 \\
H2 &    7 &     74.2 &  3.3 & 4.16e-03 &   10 &     76.2 &  3.0 & 4.87e-03 & 1.43 &   13 &     75.5 &  3.1 & 4.77e-03 & 1.86 \\
\rowcolor[gray]{0.9}
H3 &    6 &     77.4 &  1.5 & 4.08e-03 &    8 &     80.8 &  1.1 & 4.87e-03 & 1.33 &   11 &     78.9 &  1.3 & 4.48e-03 & 1.83 \\
    \hline
  \end{tabular}}

  \bigskip
  \footnotesize
  \begin{center}{ deblurring results for $\delta = 0.001$, $\varepsilon = 0.001$}\end{center}
  \hspace*{-0.7cm}{\tt \noindent\begin{tabular}{|r|rrrr|rrrrr|rrrrr|}
    \cline{2-15}
    \multicolumn{1}{c|}{} & \multicolumn{4}{c|}{IADMM ($\alpha = 0.5$)}
        & \multicolumn{5}{c|}{IADMM ($\alpha = 0.2$)} & \multicolumn{5}{c|}{ADMM(TV1)}\\
    \hline
    \rowcolor[gray]{0.9}
    IMG & I1 & ERROR & SNR & RES & I2 & ERROR & SNR & RES & I2/I1 & I3 & ERROR & SNR & RES & I3/I1 \\
L0 &   19 &     12.3 & 14.1 & 9.96e-04 &   32 &     12.2 & 14.2 & 9.75e-04 & 1.68 &   40 &     12.2 & 14.2 & 9.82e-04 & 2.11 \\
\rowcolor[gray]{0.9}
L1 &   22 &     12.0 & 12.9 & 9.92e-04 &   37 &     11.8 & 13.0 & 9.69e-04 & 1.68 &   46 &     11.8 & 13.0 & 9.82e-04 & 2.09 \\
L2 &   14 &     13.2 & 12.0 & 9.86e-04 &   23 &     13.1 & 12.1 & 1.00e-03 & 1.64 &   29 &     13.1 & 12.1 & 9.99e-04 & 2.07 \\
\rowcolor[gray]{0.9}
L3 &   11 &     12.2 & 15.9 & 9.10e-04 &   17 &     12.2 & 15.9 & 9.76e-04 & 1.55 &   21 &     12.3 & 15.8 & 9.96e-04 & 1.91 \\
M0 &   26 &     21.7 & 16.8 & 9.84e-04 &   43 &     21.5 & 16.8 & 9.78e-04 & 1.65 &   54 &     21.4 & 16.9 & 9.79e-04 & 2.08 \\
\rowcolor[gray]{0.9}
M1 &   29 &     20.3 & 16.6 & 9.91e-04 &   48 &     20.1 & 16.7 & 9.80e-04 & 1.66 &   60 &     20.1 & 16.7 & 9.86e-04 & 2.07 \\
M2 &   16 &     27.2 & 15.4 & 9.52e-04 &   26 &     27.1 & 15.4 & 9.84e-04 & 1.62 &   33 &     27.0 & 15.5 & 9.79e-04 & 2.06 \\
\rowcolor[gray]{0.9}
M3 &   24 &     25.1 & 15.5 & 9.66e-04 &   39 &     25.0 & 15.5 & 9.80e-04 & 1.62 &   48 &     25.1 & 15.5 & 1.00e-03 & 2.00 \\
H0 &   28 &     50.7 & 13.2 & 9.92e-04 &   46 &     50.3 & 13.3 & 9.89e-04 & 1.64 &   58 &     50.2 & 13.3 & 9.85e-04 & 2.07 \\
\rowcolor[gray]{0.9}
H1 &   33 &     48.5 &  9.4 & 9.82e-04 &   53 &     48.5 &  9.4 & 9.99e-04 & 1.61 &   67 &     48.3 &  9.4 & 9.91e-04 & 2.03 \\
H2 &   22 &     54.3 &  6.0 & 9.88e-04 &   36 &     53.9 &  6.0 & 9.98e-04 & 1.64 &   46 &     53.6 &  6.1 & 9.80e-04 & 2.09 \\
\rowcolor[gray]{0.9}
H3 &   18 &     57.3 &  4.1 & 9.43e-04 &   29 &     57.2 &  4.1 & 9.81e-04 & 1.61 &   36 &     57.4 &  4.1 & 9.99e-04 & 2.00 \\
    \hline
  \end{tabular}}

  \bigskip
  \footnotesize
  \begin{center}{ deblurring results for $\delta = 0.001$, $\varepsilon = 0.0005$}\end{center}
  \hspace*{-0.7cm}{\tt \noindent\begin{tabular}{|r|rrrr|rrrrr|rrrrr|}
    \cline{2-15}
    \multicolumn{1}{c|}{} & \multicolumn{4}{c|}{IADMM ($\alpha = 0.5$)}
        & \multicolumn{5}{c|}{IADMM ($\alpha = 0.2$)} & \multicolumn{5}{c|}{ADMM(TV1)}\\
    \hline
    \rowcolor[gray]{0.9}
    IMG & I1 & ERROR & SNR & RES & I2 & ERROR & SNR & RES & I2/I1 & I3 & ERROR & SNR & RES & I3/I1 \\
L0 &   36 &     10.8 & 15.2 & 4.94e-04 &   59 &     10.8 & 15.3 & 4.92e-04 & 1.64 &   73 &     10.8 & 15.3 & 4.99e-04 & 2.03 \\
\rowcolor[gray]{0.9}
L1 &   42 &     10.5 & 14.0 & 4.96e-04 &   68 &     10.5 & 14.1 & 4.99e-04 & 1.62 &   86 &     10.5 & 14.1 & 4.95e-04 & 2.05 \\
L2 &   27 &     11.5 & 13.2 & 4.91e-04 &   44 &     11.4 & 13.3 & 4.95e-04 & 1.63 &   55 &     11.4 & 13.3 & 4.97e-04 & 2.04 \\
\rowcolor[gray]{0.9}
L3 &   21 &     10.6 & 17.2 & 4.80e-04 &   34 &     10.5 & 17.2 & 4.88e-04 & 1.62 &   42 &     10.6 & 17.1 & 4.97e-04 & 2.00 \\
M0 &   49 &     18.2 & 18.3 & 4.99e-04 &   80 &     18.1 & 18.3 & 4.96e-04 & 1.63 &  100 &     18.1 & 18.3 & 4.97e-04 & 2.04 \\
\rowcolor[gray]{0.9}
M1 &   50 &     17.4 & 17.9 & 5.50e-04 &   79 &     17.5 & 17.9 & 5.64e-04 & 1.58 &  113 &     16.8 & 18.2 & 4.96e-04 & 2.26 \\
M2 &   34 &     22.1 & 17.2 & 5.83e-04 &   54 &     22.1 & 17.2 & 5.90e-04 & 1.59 &  158 &     17.8 & 19.1 & 5.00e-04 & 4.65 \\
\rowcolor[gray]{0.9}
M3 &   41 &     22.0 & 16.6 & 6.65e-04 &   65 &     22.0 & 16.6 & 6.74e-04 & 1.59 &  401 &     16.3 & 19.2 & 5.00e-04 & 9.78 \\
H0 &   40 &     45.8 & 14.1 & 6.73e-04 &   64 &     45.8 & 14.1 & 6.85e-04 & 1.60 &  110 &     41.8 & 14.9 & 4.96e-04 & 2.75 \\
\rowcolor[gray]{0.9}
H1 &   68 &     39.5 & 11.1 & 4.99e-04 &  108 &     39.6 & 11.1 & 5.05e-04 & 1.59 &  138 &     39.3 & 11.2 & 4.97e-04 & 2.03 \\
H2 &   42 &     45.3 &  7.6 & 4.90e-04 &   67 &     45.3 &  7.6 & 4.99e-04 & 1.60 &   84 &     45.3 &  7.6 & 4.99e-04 & 2.00 \\
\rowcolor[gray]{0.9}
H3 &   32 &     48.8 &  5.5 & 4.91e-04 &   52 &     48.6 &  5.5 & 4.94e-04 & 1.62 &   65 &     48.7 &  5.5 & 4.96e-04 & 2.03 \\
    \hline
  \end{tabular}}
    \caption{Deblurring results for a small value of $\delta$  ($\delta = 0.001$) using different values of the tolerance $\varepsilon$. \emph{Iteration number} (\texttt{I1, I2, I3}), \emph{Real error} (\texttt{ERROR}), \emph{SNR}, \emph{Residual} (\texttt{RES}) and \emph{efficiency rates} \texttt{I2/I1} and \texttt{I3/I1} using the IADMM ($\alpha=0.5$ and $0.2$) and the ADMM (TV1) methods applied to all the 12 test images (\texttt{IMG}). \label{table1}}
  \end{table*}

The proposed algorithm IADMM (Algorithm~\ref{algo}) is compared with the widely used augmented Lagrangian
methods (ADMM \cite{chen2015general}) for image deblurring. We mainly consider two models, i.e., $\varphi(\cdot)=|\cdot|$ and $\varphi(\cdot)=|\cdot|^{\frac{1}{2}}$. We call them as TV1 and TV(1/2) methods, respectively. Note that TV1 is a convex method, while TV(1/2) is nonconvex. In the tests we have considered (unless so indicated) for all the methods the value of the penalty parameter $\delta = 0.001$ (a small one) and/or $\delta = 10$ (a large one) just to see the behaviour of the IADMM algorithm.

The
performance of the deblurring algorithms is quantitatively measured by means of the \emph{objective function} (Equations \ref{debblurring} or \ref{debblurring2}), the \emph{Real Error} as $\| \cdot \|_2$ of the difference between the original and deblurred images, the \emph{signal-to-noise ratio (SNR)} \cite{chen2015general}
\begin{align}
\textrm{SNR}(u,u^*)=10 \times \log_{10}(\frac{\|u-\bar{u}\|^2}{\|u-u^*\|^2}),
\end{align}
where  $u$ and $u^*$ denote the original image and the restored image, respectively, and $\bar{u}$ represents the mean of the original image $u$, and the \emph{residual}
($\texttt{res}(k, k+1)$ and ${\texttt{res}_i}(k, k+1)$ in the standard and inertial versions, respectively) as described in \cite{chen2015general}:
\begin{equation}
\begin{array}{l}
\texttt{res}(k, k+1) := \displaystyle\frac{\| (u^{k+1}, \, p^{k+1})-(u^{k}, \, p^{k})\|}{1+\|(u^{k}, \, p^{k})\|}, \\[3.ex]
\texttt{res}_{i}(k, k+1) := \displaystyle\frac{\| (u^{k+1}, \, p^{k+1})-(\widehat{u}^{k}, \, \widehat{p}^{k})\|}{1+\|(\widehat{u}^{k}, \, \widehat{p}^{k})\|}.
\end{array}
\label{residu}
\end{equation}
In the tests we do not provide CPU tests as all the algorithms have a very similar computation cost per iteration (mainly from the FFT routines). Therefore, there is almost no difference between pictures showing iterations or CPU cost, and the consumed CPU is basically proportional to the respective number of iterations.

On our first test we use the brain M0 image with a low value of $\delta = 0.001$ and we show in Figure \ref{tests} that the performances of the ADMM TV1 and TV(1/2) are quite similar in error and SNR. Therefore, in the rest of comparisons we will just consider the TV1 method. For the IADMM method the inertial parameter $\alpha$ in Algorithm~\ref{algo} is investigated firstly by using two different values ($\alpha = 0.2$ and $0.5$).
The main difference observed in these tests is that the nonconvex method is more unstable once reached the maximum precision (at this point the ADMM TV1 convex method seems to be the more stable with a quite smooth behaviour). The fastest convergence is observed using the IADMM (with the largest \emph{stepsize} value $\alpha = 0.5$) method, but when the maximum precision is obtained unstable behaviour appears. Therefore, using mainly the information provided by the \texttt{residual} (Eq. (\ref{residu})), we provide a \emph{stop control} criterion (in the same spirit as
\cite{chen2015general})
that stops the iterative process at the black points of Figure \ref{tests} (d) in the tests. That is, we stop when the following condition is hold
\begin{equation}
\label{stopc}
\texttt{res}(k, k+1) < \varepsilon \quad \textrm{or} \quad \texttt{res}(k-1, k) < \texttt{res}(k, k+1).
\end{equation}
In the first case the convergence till the desired tolerance error is obtained, while in the second case the algorithm has reached its stability limit and the residual grows, behaving later in an unstable way. Note that the residual is used in the stop criterion, as it uses known data from the iterations (and which does not depend on the original image that is unknown). We remark that the use of stop control techniques avoids the use of unnecessary iterations, and also to stop at the limit accuracy of the method. Also, from the pictures we observe that the IADMM algorithm provides
enough precision in a lower number of iterations. The larger $\alpha$ means faster method, but at the price of a more unstable method as it can be seen on Fig. \ref{tests}(d). On that picture we observe that when the residual begins to behave chaotically, with sudden increases, it means that it is advisable to stop the iterative process as considered in the stop criterion (\ref{stopc}) (black dot points on Fig. \ref{tests}(d)). With that criterion, the IADMM method seems to be an interesting option for fast deblurring problems.

\begin{figure}[htb!]
\centering
 \includegraphics[width=0.9\columnwidth]{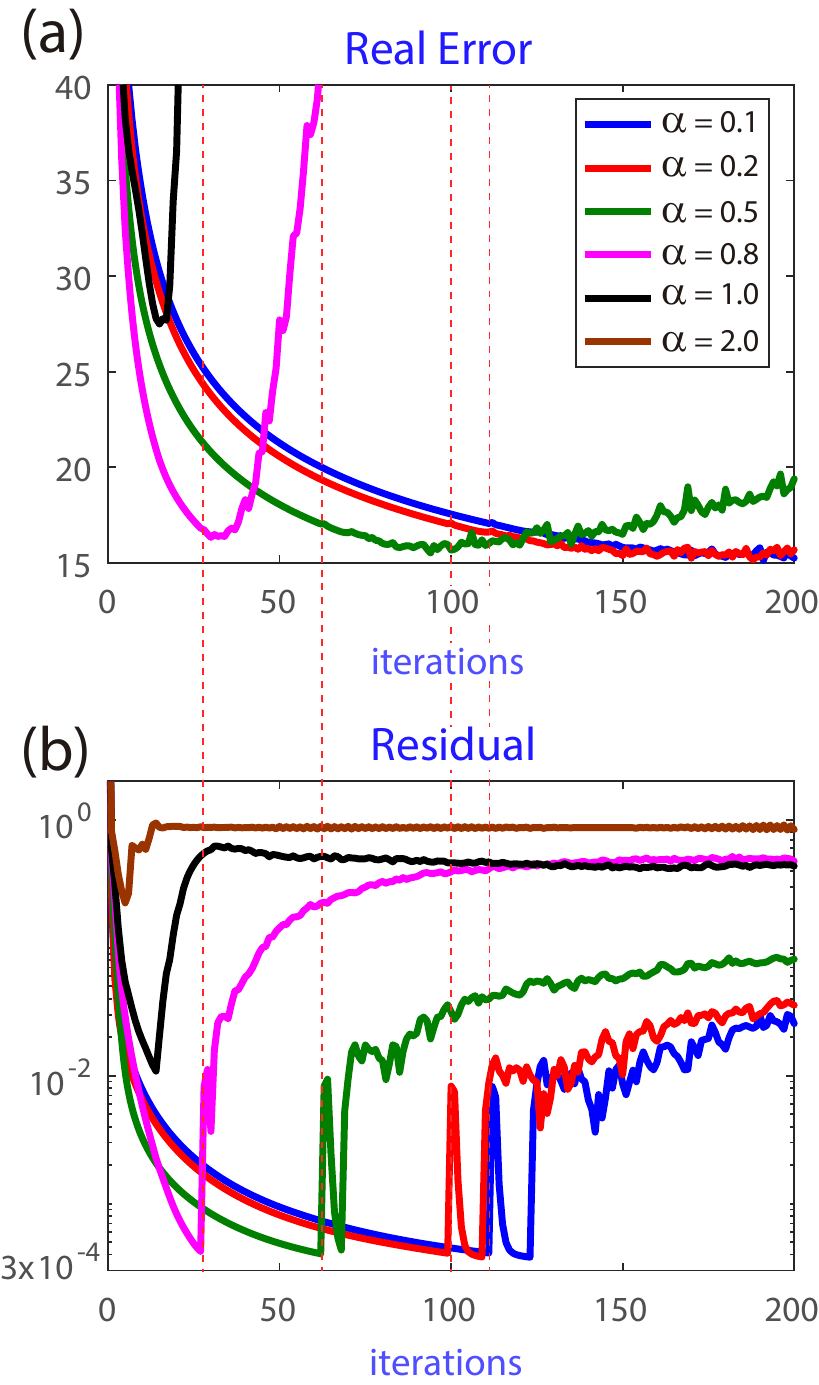}
  \caption{The evolution curves of (a) the \emph{Real Error} and (b) the \emph{Residual}, for the brain M0 image versus the iteration number. All the simulations have been done using the IADMM algorithm with several values of the parameter $\alpha = 0.1, 0.2, 0.5, 0.8, 1.0$ and $2.0$ and $\delta = 0.001$.}
  \label{testsalpha}
\end{figure}

\begin{figure}[htb!]
\centering
 \includegraphics[width=0.9\columnwidth]{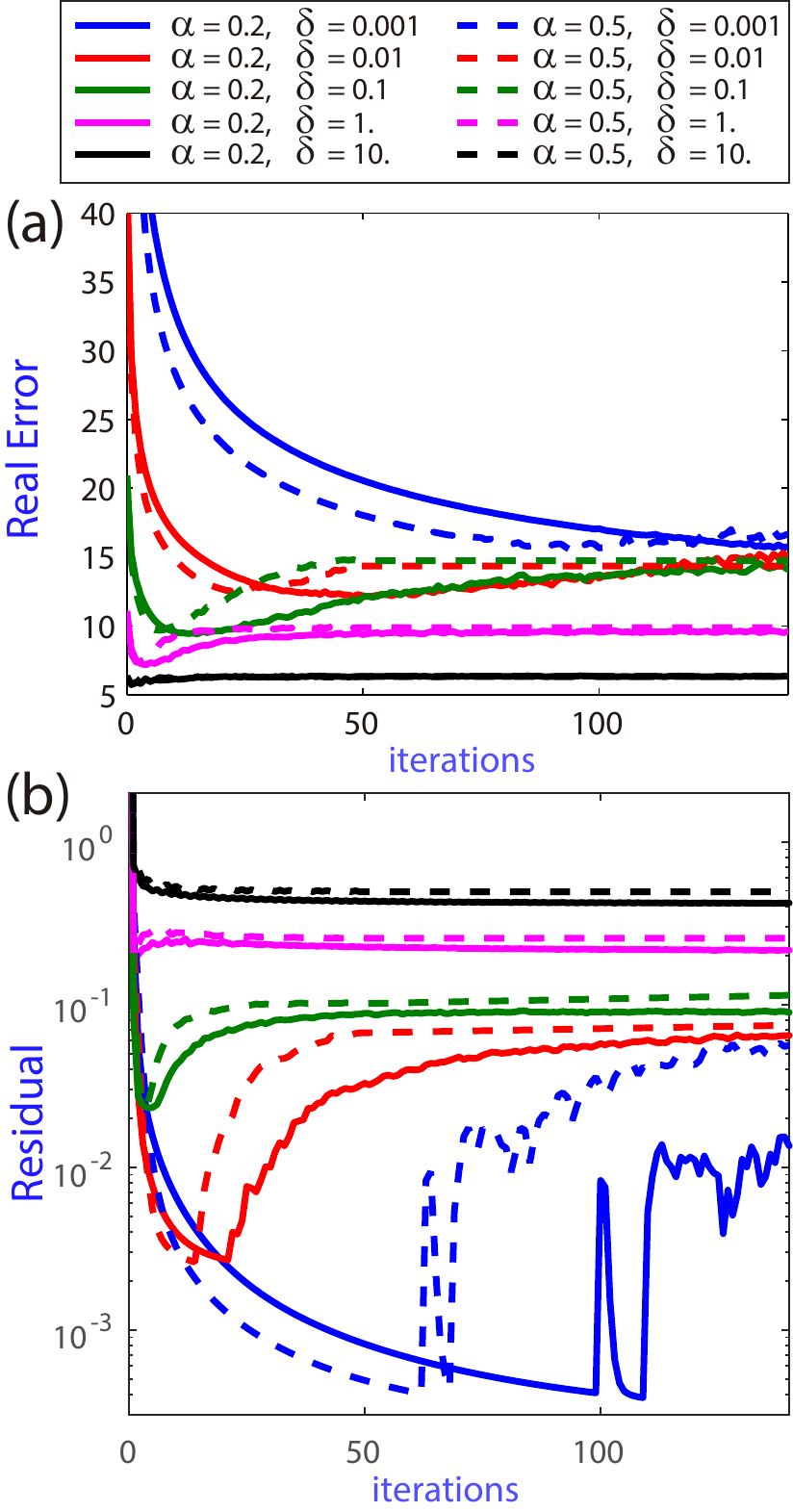}
  \caption{The evolution curves of (a) the \emph{Real Error} and (b) the \emph{Residual}, for the brain M0 image versus the iteration number. All the simulations have been done using the IADMM algorithm with several values of the parameters $\alpha = 0.2, 0.5$ and $\delta = 0.001, 0.01, 0.1, 1$ and $10$.}
  \label{testsdelta}
\end{figure}

To observe more clearly the influence of the parameter $\alpha$ in Algorithm~\ref{algo} we perform several tests on the brain M0 image on Figure \ref{testsalpha} for values $\alpha = 0.1, 0.2, 0.5, 0.8, 1.0$ and $2.0$. Note that this parameter plays a role similar to the \emph{stepsize} (as it also occurs to the parameter $\delta$), as it controls the perturbation at each step. A large value will provide, when the method works, a quite fast method, but on the other hand it makes the method more unstable. In fact, from the plot \ref{testsalpha}(b) we observe that in this case it would be optimal to use the parameter value $\alpha = 0.8$ in combination with the stop criterion, giving the maximum precision in just 28 iterations. Besides, it is shown that after the values selected by the stop criterion the residual begins to oscillate among values that provides similar error but that generates an unstable behaviour giving rise to an increment of the error in subsequent iterations (this instability is delayed when the parameter $\alpha$ decreases, what is expected because the increment is smaller, as the vertical lines connecting the error and residual plots show).

The influence of the penalty parameter $\delta$ is also quite relevant, but a detailed analysis is out of the scope of this paper. On the Figure~\ref{testsdelta} we show the evolution of the residual using two values of $\alpha = 0.2, 0.5$
and several values of the parameter $\delta = 0.001, 0.01, 0.1, 1$ and $10$. We observe that low values of the penalty parameter $\delta$ gives a lower residual, but the error is lower for large $\delta$ providing a faster convergence, and it has a big effect on the empirical performance of the methods as shown in \cite{xu2016empirical}, but it remains to study optimal combinations of the parameters $\delta$ and $\alpha$ and suitable criteria for an automatic selection (this will be part of the next steps in our study of these methods).

\begin{figure}[htb!]
\centering
 \includegraphics[width=0.9\columnwidth]{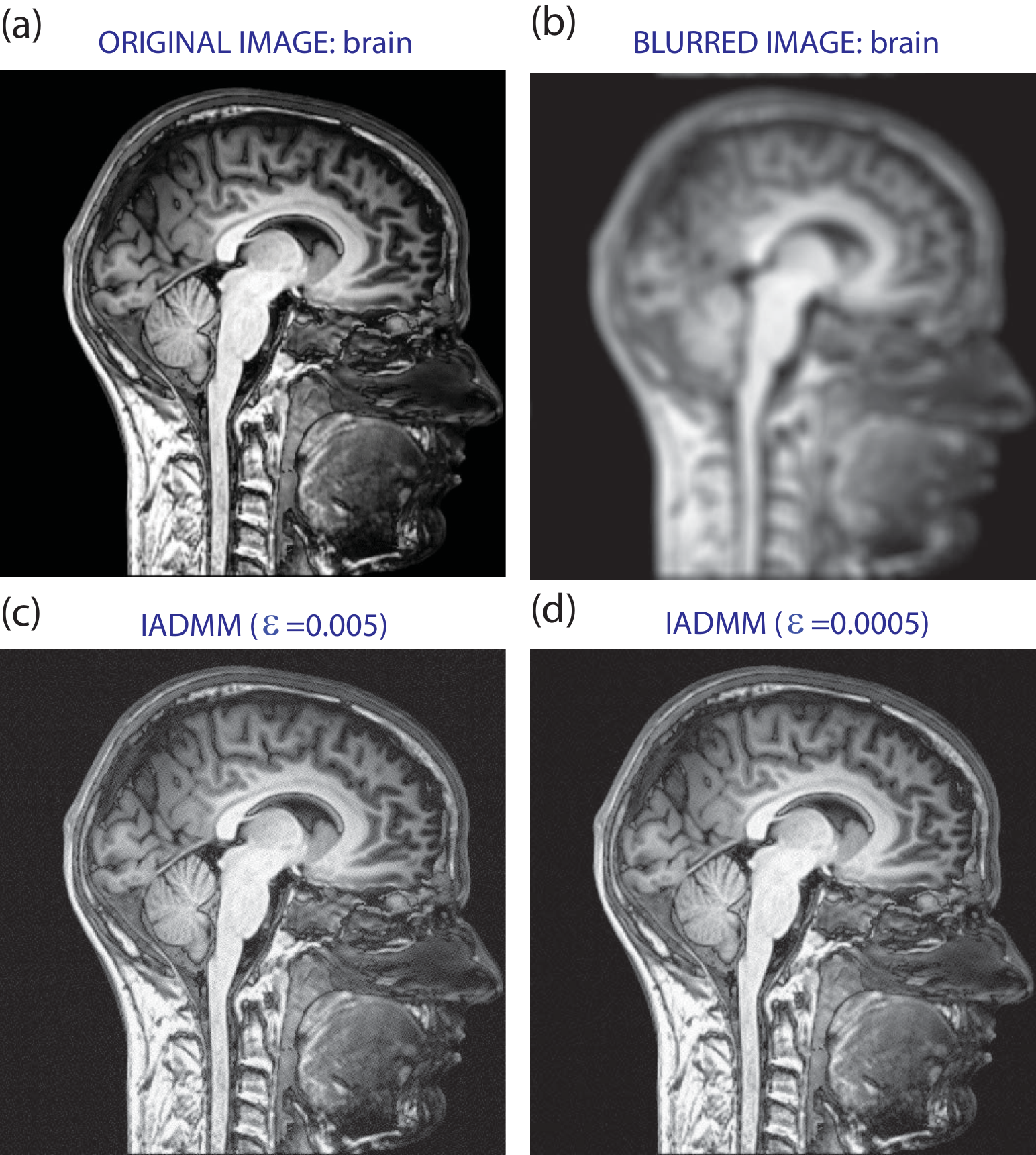}
  \caption{Deblurred images at different stages of the IADMM method for brain M0 image. (a) Original image; (b) Blurred image;  (c) Recovery by IADMM using error tolerance $\varepsilon = 5 \times 10^{-3}$; (d) Recovery by IADMM using error tolerance $\varepsilon = 5 \times 10^{-4}$. \label{fig4brains}}
\end{figure}

On Figure \ref{fig4brains}, we present the original medium resolution brain M0 image, the blurred one using, as indicated, a convolution with Gaussian kernel, and the results of the IADMM deblurring images using two error tolerances ($\varepsilon=5 \times 10^3$, $5 \times 10^4$) in the stop criterion (\ref{stopc}). We can see that in both cases the quality of the recovered image is visually good.

\begin{table*}[htb!]
  \footnotesize
  \begin{center}{deblurring results for $\delta = 10$, $\varepsilon = 0.01$}\end{center}
  \hspace*{-0.7cm}{\tt \noindent\begin{tabular}{|r|rrrr|rrrrr|rrrrr|}
    \cline{2-15}
    \multicolumn{1}{c|}{} & \multicolumn{4}{c|}{IADMM ($\alpha = 0.5$)}
        & \multicolumn{5}{c|}{IADMM ($\alpha = 0.2$)} & \multicolumn{5}{c|}{ADMM(TV1)}\\
    \hline
    \rowcolor[gray]{0.9}
    IMG & I1 & ERROR & SNR & RES & I2 & ERROR & SNR & RES & I2/I1 & I3 & ERROR & SNR & RES & I3/I1 \\
L0 &    5 &      3.0 & 26.5 & 4.72e-01 &    5 &      2.8 & 27.1 & 4.31e-01 & 1.00 &   10 &      6.2 & 20.0 & 9.56e-02 & 2.00 \\
\rowcolor[gray]{0.9}
L1 &    3 &      3.3 & 24.1 & 5.38e-01 &    5 &      3.2 & 24.3 & 4.62e-01 & 1.67 &   10 &      5.8 & 19.2 & 9.48e-02 & 3.33 \\
L2 &    1 &      3.8 & 22.9 & 8.22e-01 &    3 &      3.7 & 23.0 & 5.13e-01 & 3.00 &   10 &      9.3 & 15.0 & 9.68e-02 & 10.00 \\
\rowcolor[gray]{0.9}
L3 &    6 &      1.4 & 34.6 & 2.41e-01 &    7 &      1.3 & 35.5 & 2.00e-01 & 1.17 &   10 &      7.0 & 20.7 & 9.86e-02 & 1.67 \\
M0 &    3 &      6.2 & 27.6 & 5.79e-01 &    5 &      6.0 & 28.0 & 4.98e-01 & 1.67 &   10 &      7.0 & 26.5 & 9.74e-02 & 3.33 \\
\rowcolor[gray]{0.9}
M1 &    7 &      4.3 & 30.2 & 4.40e-01 &    7 &      3.9 & 31.0 & 4.12e-01 & 1.00 &   10 &      5.4 & 28.1 & 9.78e-02 & 1.43 \\
M2 &    4 &      3.0 & 34.5 & 1.50e-01 &    4 &      3.0 & 34.7 & 1.00e-01 & 1.00 &   10 &     20.6 & 17.8 & 9.90e-02 & 2.50 \\
\rowcolor[gray]{0.9}
M3 &    1 &     14.6 & 20.2 & 8.94e-01 &    1 &     14.5 & 20.2 & 8.68e-01 & 1.00 &   10 &     31.7 & 13.4 & 9.64e-02 & 10.00 \\
H0 &    5 &     13.8 & 24.5 & 4.96e-01 &    5 &     12.6 & 25.3 & 4.44e-01 & 1.00 &   10 &     15.2 & 23.7 & 9.52e-02 & 2.00 \\
\rowcolor[gray]{0.9}
H1 &    3 &     16.6 & 18.7 & 5.73e-01 &    5 &     16.2 & 18.9 & 4.98e-01 & 1.67 &   10 &     20.4 & 16.9 & 9.47e-02 & 3.33 \\
H2 &    5 &     15.7 & 16.8 & 5.11e-01 &    5 &     14.6 & 17.4 & 4.57e-01 & 1.00 &   10 &     20.7 & 14.4 & 9.29e-02 & 2.00 \\
    \hline
  \end{tabular}}
    \caption{Deblurring results for a large value of $\delta$  ($\delta = 10$)  using the tolerance $\varepsilon = 0.01$. \emph{Iteration number} (\texttt{I1, I2, I3}), \emph{Real error} (\texttt{ERROR}), \emph{SNR}, \emph{Residual} (\texttt{RES}) and \emph{efficiency rates} \texttt{I2/I1} and \texttt{I3/I1} using the IADMM ($\alpha=0.5$ and $0.2$) and the ADMM (TV1) methods applied to all the 12 test images (\texttt{IMG}). \label{table2}}
  \end{table*}

  \begin{table*}[htb!]
  \footnotesize
  \begin{center}{deblurring results for $\varepsilon = 0.01$ on M0 image}\end{center}
  \hspace*{-0.7cm}{\tt \noindent\begin{tabular}{|r|rrrr|rrrrr|rrrrr|}
    \cline{2-15}
    \multicolumn{1}{c|}{} & \multicolumn{4}{c|}{IADMM ($\alpha = 0.5$)}
        & \multicolumn{5}{c|}{IADMM ($\alpha = 0.2$)} & \multicolumn{5}{c|}{ADMM(TV1)}\\
    \hline
    \rowcolor[gray]{0.9}
    $\delta$ & I1 & ERROR & SNR & RES & I2 & ERROR & SNR & RES & I2/I1 & I3 & ERROR & SNR & RES & I3/I1 \\
0.001 &    6 &     33.5 & 13.0 & 7.40e-03 &    8 &     35.3 & 12.5 & 8.70e-03 & 1.33 &    9 &     34.0 & 12.8 & 7.50e-03 & 1.50 \\
\rowcolor[gray]{0.9}
0.01 &    4 &     19.1 & 17.9 & 8.88e-03 &    4 &     21.2 & 16.9 & 9.82e-03 & 1.00 &    5 &     20.0 & 17.4 & 7.68e-03 & 1.25 \\
0.10 &    3 &     11.2 & 22.5 & 2.39e-02 &    4 &     11.4 & 22.3 & 2.31e-02 & 1.33 &    5 &     10.8 & 22.8 & 2.32e-02 & 1.67 \\
\rowcolor[gray]{0.9}
1.00 &    1 &      8.0 & 25.4 & 2.25e-01 &    2 &      7.5 & 26.0 & 2.03e-01 & 2.00 &    3 &      7.3 & 26.3 & 2.20e-01 & 3.00 \\
10.00 &    3 &      6.2 & 27.6 & 5.79e-01 &    5 &      6.0 & 28.0 & 4.98e-01 & 1.67 &   10 &      7.0 & 26.5 & 9.74e-02 & 3.33 \\
    \hline
  \end{tabular}}
    \caption{Deblurring results for several values of parameter $\delta$ using the tolerance $\varepsilon = 0.01$ on M0 image. \emph{Iteration number} (\texttt{I1, I2, I3}), \emph{Real error} (\texttt{ERROR}), \emph{SNR}, \emph{Residual} (\texttt{RES}) and \emph{efficiency rates} \texttt{I2/I1} and \texttt{I3/I1} using the IADMM ($\alpha=0.5$ and $0.2$) and the ADMM (TV1) methods. \label{table3}}
  \end{table*}

A more detailed analysis is shown on Table \ref{table1}, where we present for the IADMM ($\alpha=0.5$ and $0.2$) and the ADMM (TV1) methods (we do not show results for the TV(1/2) as they are quite similar to those of the TV1 ones), the iteration number, real error (\texttt{ERROR}), \texttt{SNR}, the residual (\texttt{RES}) and the \emph{efficiency rates} of the IADMM ($\alpha=0.5$) method over the two other ones, \texttt{I2/I1} and \texttt{I3/I1}, for different values of the tolerance $\varepsilon$ and using a small value of $\delta=0.001$. In all the methods we have used the same values of the parameters and a similar stop criterion (\ref{stopc}) with the respective definition of the residual (Eqn. (\ref{residu})). From the tests, the inertial IADMM nonconvex methods are the fastest, as expected, but also with a more unstable behaviour, also as expected due to the nonconvex version of the methods. Note that from Figure~\ref{tests} the ADMM TV1 and TV(1/2) perform similarly for low $\delta$, but if we make simulations small differences appear in medium-high resolution, obtaining for the M1 image the TV1 method 113 iterations and a ratio 2.26, but the TV(1/2)  99 iterations and a ratio  1.98; and for the H0 image the TV1 uses 110 iterations and a ratio  2.75, but the TV(1/2)   81 iterations and a ratio  2.03. That is, the nonconvex TV(1/2) may perform better in some circumstances, but the instability may also appear. In case of using larger values of the parameter $\delta$, as the convergence theorems suggest, all the methods are much faster, giving a
better performance as shown on Table~\ref{table2} for $\delta = 10$. Finally, we present on Table~\ref{table3} some results obtained by changing the value of the parameter $\delta$ and with the fixed value of the tolerance $\varepsilon = 0.01$. Again, the IADMM ($\alpha=0.5$) method presents the best performance. Therefore, globally, the inertial IADMM version seems to be an interesting option, being the fastest one.

As above commented, how to select optimal combinations of the parameters $\delta$ and $\alpha$, the use on other types of blur, development of suitable criteria for automatic selection of all the parameters and so on, remains a research issue, but these goals are out of the scope of this article.

\section{Conclusion}
\label{sec:6}
In this paper, a more efficient  nonconvex inertial alternating minimization algorithm (IADMM)
is developed for solving the Total Variation model for  image deblurring. The proposed scheme is based on using the inertial proximal strategy on nonconvex ADMM methods. By the using the Kurdyka-{\L}ojasiewicz property, we prove the convergence of the algorithm under several reasonable assumptions.
 Numerical experiments demonstrate that
the proposed algorithm overall outperforms the widely used augmented Lagrangian ADMM methods, being a fast option, although it can be unstable for high precision. This instability is easily controlled via a suitable stop control criteria.

\section*{Appendix A: proof of Lemma \ref{lemmaini}}

Direct computation yields
\footnotesize
\begin{eqnarray}\label{lemmaini-t1}
\mathcal{T}^* \mathcal{T}=\left(
                  \begin{array}{cc}
                    (I_N\otimes D)(I_N\otimes D)^* &  \\
                     & (D\otimes I_N)(D\otimes I_N)^* \\
                  \end{array}
                \right).
\end{eqnarray}
\normalsize
 Then, we have
 \footnotesize
\begin{align}\label{lemmaini-t2}
   &\lambda_{\min}( \mathcal{T}^*)^2=\lambda_{\min}( \mathcal{T}\mathcal{T}^*)=\\
   &\min\big\{\lambda_{\min}\big( (I_N\otimes D)(I_N\otimes D)^*\big), \lambda_{\min}\big((D\otimes I_N)(D\otimes I_N)^* \big)\big\}. \nonumber
\end{align}
\normalsize
With Lemma 7.2 in \cite{alizadeh1998primal}, we obtain
\footnotesize
\begin{align}\label{lemmaini-t2+}
  &\lambda_{\min}\big( (I_N\otimes D)(I_N\otimes D)^*\big)\nonumber\\
  &\quad\quad=\lambda_{\min}\big((D\otimes I_N)(D\otimes I_N)^* \big)=\lambda_{\min}(DD^*).
\end{align}
\normalsize
Noting $DD^*$ is a symmetric tridiagonal matrix
\footnotesize
\begin{equation}\label{DD}
    DD^*=\left(
           \begin{array}{cccc}
             2 & -1 &  &  \\
             -1 & 2 & \ddots &  \\
              & \ddots & \ddots & -1 \\
              &  & -1 & 2 \\
           \end{array}
         \right),
\end{equation}
\normalsize
and using the result in \cite{link},
\footnotesize
\begin{align}
    &\lambda_{\min}(DD^*)=2+2\cos\bigg(\pi-\frac{\pi}{N}\bigg)\nonumber\\
    &\quad\quad=4\cos^2\bigg(\frac{\pi}{2}-\frac{\pi}{2N}\bigg)=4\sin^2\bigg(\frac{\pi}{2N}\bigg),
\end{align}
\normalsize
we obtain
$
    \lambda_{\min}(\mathcal{T}^*)\geq 2\sin\bigg(\frac{\pi}{2N}\bigg).
$
Now, taking into account that $\mathcal{T}^*$ is full row-rank,
$
  \|\mathcal{T}^* x\|\geq \lambda_{\min}(\mathcal{T}^*)\|x\|\geq 2\sin\bigg(\frac{\pi}{2N}\bigg)\|x\|.
$
%%%%%%%%%%%%%%%%%%%%%%%%%%%%%%%%%%%%%%%%%%%%%%%%
\section*{Appendix B: proof of Lemma \ref{fp}}
As the point $u^{k+1}$ is the minimization of $\mathcal{L}_{\delta}^{\varphi}(u,v^{k+1},\widehat{p}^k)$, then the optimization condition yields
\footnotesize
\begin{equation}
    K^*(K u^{k+1}-f)-\mathcal{T}^{*}\widehat{p}^{k}+\delta \mathcal{T}^{*}(\mathcal{T} u^{k+1}-v^{k+1})=\textbf{0}.
\end{equation}
\normalsize
Therefore, we can derive
\footnotesize
\begin{equation}\label{fp+t2}
    K^*(K u^{k+1}-f)-\mathcal{T}^{*}p^{k+1}=\textbf{0},
\end{equation}
\normalsize
and replacing $k+1$ by $k$
\footnotesize
\begin{equation}\label{boundused}
 K^*(K u^{k}-f)-\mathcal{T}^{*}p^{k}=\textbf{0}.
\end{equation}
\normalsize
Subtracting (\ref{fp+t2}) and (\ref{boundused}), and using Lemma~\ref{lemmaini},
\footnotesize
\begin{align}
&\|p^k-p^{k+1}\|\leq \theta\|\mathcal{T}^{*}(p^k-p^{k+1})\|\nonumber\\
&\quad=\theta\|K^*K(u^{k+1}-u^{k}\|\leq\theta\|K\|_2^2\cdot\|u^{k+1}-u^{k}\|.
\end{align}
\normalsize
%%%%%%%%%%%%%%%%%%%%%%%%%%%%%%%%%%%%%%%%%%%%%%%%%%
\section*{Appendix C: proof of Lemma \ref{descend1}}
Note that $v^{k+1}$ is the minimizer of $\mathcal{L}_{\delta}^{\varphi}(u^k,v,\widehat{p}^k)$ with respect to the variable $v$. Then we have:
\footnotesize
\begin{eqnarray}\label{s1temp1}
\mathcal{L}_{\delta}^{\varphi}(u^k,v^{k+1},\widehat{p}^k)\leq \mathcal{L}_{\delta}^{\varphi}(u^k,v^k,\widehat{p}^k).
\end{eqnarray}
\normalsize
In (\ref{des+}), we have obtained
\footnotesize
\begin{align}\label{s1temp2}
\mathcal{L}_{\delta}^{\varphi}(u^{k+1},v^{k+1},\widehat{p}^k)&+\frac{\nu}{2}\|u^{k+1}-u^k\|^2\leq \mathcal{L}_{\delta}^{\varphi}(u^k,v^{k+1},\widehat{p}^k).
\end{align}
\normalsize
By direct calculations, we obtain
\footnotesize
\begin{align}\label{s1temp3}
&\mathcal{L}_{\delta}^{\varphi}(u^{k+1},v^{k+1},p^{k+1})\nonumber\\
&=\mathcal{L}_{\delta}^{\varphi}(u^{k+1},v^{k+1},\widehat{p}^k)+\langle \widehat{p}^k-p^{k+1},\mathcal{T} u^{k+1}-v^{k+1}\rangle\nonumber\\
    &=\mathcal{L}_{\delta}^{\varphi}(u^{k+1},v^{k+1},\widehat{p}^k)+\frac{1}{\delta}\|\widehat{p}^k-p^{k+1}\|_2^2\nonumber\\
   &\leq\mathcal{L}_{\delta}^{\varphi}(u^{k+1},v^{k+1},\widehat{p}^k)+\frac{2}{\delta}\|\widehat{p}^k-p^{k}\|_2^2+\frac{2}{\delta}\|p^{k}-p^{k+1}\|_2^2\nonumber\\
   &\leq\mathcal{L}_{\delta}^{\varphi}(u^{k+1},v^{k+1},\widehat{p}^k)+\frac{2\alpha^2}{\delta}\|p^{k}-p^{k-1}\|_2^2+\frac{2}{\delta}\|p^{k}-p^{k+1}\|_2^2\nonumber\\
   &\leq\mathcal{L}_{\delta}^{\varphi}(u^{k+1},v^{k+1},\widehat{p}^k)+\frac{2\theta^2\|K\|_2^4}{\delta}\|u^{k+1}-u^{k}\|^2\nonumber\\
   & \qquad \qquad + \frac{2\alpha^2\theta^2\|K\|_2^4}{\delta}\|u^{k}-u^{k-1}\|^2.
\end{align}
\normalsize
and
\footnotesize
\begin{align}\label{s1temp4}
&\mathcal{L}_{\delta}^{\varphi}(u^{k},v^{k},\widehat{p}^{k})=\mathcal{L}_{\delta}^{\varphi}(u^{k},v^{k},p^k)+\langle p^k-\widehat{p}^k,\mathcal{T} u^{k}-v^{k}\rangle\nonumber\\
    &=\mathcal{L}_{\delta}^{\varphi}(u^{k},v^{k},p^k)+\frac{1}{\delta}\langle p^k-\widehat{p}^k,p^{k+1}-\widehat{p}^k\rangle\nonumber\\
    &\leq\mathcal{L}_{\delta}^{\varphi}(u^{k},v^{k},p^k)+\frac{\|p^k-\widehat{p}^k\|^2}{2\delta}+\frac{\|p^{k+1}-\widehat{p}^k\|^2}{2\delta}\nonumber\\
    &\leq\mathcal{L}_{\delta}^{\varphi}(u^{k},v^{k},p^k)+\frac{3\alpha^2\|p^k-p^{k-1}\|^2}{2\delta}+\frac{\|p^{k+1}-p^k\|^2}{\delta}\nonumber\\
    &\leq\mathcal{L}_{\delta}^{\varphi}(u^{k},v^{k},p^k)+ \frac{\theta^2\|K\|_2^4}{\delta}\|u^{k+1}-u^{k}\|^2\nonumber\\
    &\qquad \qquad +\frac{3\alpha^2\theta^2\|K\|_2^4}{2\delta}\|u^{k}-u^{k-1}\|^2.
\end{align}
\normalsize
Combining the above equations, we obtain
\footnotesize
\begin{align}\label{s1temp5}
    &\mathcal{L}_{\delta}^{\varphi}(x^{k+1},y^{k+1},p^{k+1})+\frac{\nu}{2}\|u^{k+1}-u^k\|^2\nonumber\\
    & \quad \leq\mathcal{L}_{\delta}^{\varphi}(x^{k},y^{k},p^{k})+ \frac{3\theta^2\|K\|_2^4}{\delta}\|u^{k+1}-u^{k}\|^2 \nonumber\\
    &\qquad \qquad+\frac{7\alpha^2\theta^2\|K\|_2^4}{2\delta}\|u^{k}-u^{k-1}\|^2.\nonumber
\end{align}
\normalsize
Thus, taking into account the definition of $F$ (\ref{fun}), we have
\footnotesize
\begin{eqnarray}\label{s1temp6}
    &&F(w^{k})-F(w^{k+1})\\
    &&\qquad  \geq\bigg[\frac{\nu}{2}-\bigg( \frac{3\theta^2\|K\|_2^4}{\delta}+\frac{7\alpha^2\theta^2\|K\|_2^4}{2\delta}\bigg)\bigg]\|u^{k+1}-u^k\|^2. \nonumber
\end{eqnarray}
\normalsize
Now, using the condition (\ref{delta1}) and the definition (\ref{auxh}) of $\widehat{h}$, we obtain the descent result given in Eq.~(\ref{descentcond}).

The boundness of $\{(u^k,v^k,p^k)\}_{k=0,1,2,\ldots}$ implies the boundness of $\{w^k\}_{k=0,1,2,\ldots}$, and the continuity of $F$ implies that $\{F(w^k)\}_{k=0,1,2,\ldots}$ is bounded. From the above proven result Eq.~(\ref{descentcond}),  $F(w^k)$ is decreasing. Thus, the sequence $\{F(w^k)\}_{k=0,1,2,\ldots}$ is convergent, i.e., $\lim_{k}[F(w^k)-F(w^{k+1})]=0$.  Therefore, from (\ref{descentcond}), we have
\footnotesize
\begin{equation}
    \lim_{k}\|u^{k+1}-u^k\|_2\leq\lim_{k}\sqrt{\frac{F(w^k)-F(w^{k+1})}{\widehat{h}}}=0,
\end{equation}
\normalsize
which indicates  $\lim_k \|p^{k+1}-p^k\|=0$ due to (\ref{condition1}).  The scheme of updating $p^{k+1}$ gives us $\lim_k \|v^{k+1}-v^k\|=0$.

%%%%%%%%%%%%%%%%%%%%%%%%%%%%%%%%%%%%%%%
\section*{Appendix D: proof of Lemma \ref{bound}}
From (\ref{boundused}), we derive
\footnotesize
\begin{align}\label{boundt1}
\|p^{k}\|^2\leq \theta^2\|K\|_2^4\cdot\|Ku^{k}-f\|^2.
\end{align}
\normalsize
By direct calculations, we obtain
\footnotesize
\begin{align}
&\mathcal{L}_{\delta}^{\varphi}(u^{k},v^{k},p^k)= \frac{1}{2}\|Ku^{k}-f\|^2+\sigma\|v^{k}\|_{\varphi}\nonumber\\
&\qquad -\langle p^{k},\mathcal{T} u^{k}-v^{k}\rangle+\frac{\delta}{2}\|\mathcal{T} u^{k}-v^{k}\|^2\nonumber\\
&=\frac{1}{2}\|Ku^{k}-f\|^2+\sigma\|v^{k}\|_{\varphi} +\frac{\delta}{2}\bigg\|\mathcal{T} u^{k}-v^{k}-\frac{p^{k}}{\delta}\bigg\|^2-\frac{\|p^{k}\|^2}{2\delta}\nonumber\\
&\geq \frac{1}{2}\|Ku^{k}-f\|^2+\sigma\|v^{k}\|_{\varphi}+\frac{\delta}{2}\bigg\|\mathcal{T} u^{k}-v^{k}-\frac{p^{k}}{\delta}\bigg\|^2 \nonumber\\ &\qquad -\frac{\theta^2\|K\|_2^4}{2\delta}\|Ku^{k}-f\|^2\nonumber\\
&=\bigg(\frac{1}{2}-\frac{\theta^2\|K\|_2^4}{2\delta}\bigg)\|Ku^{k}-f\|^2 +\frac{\delta}{2}\bigg\|\mathcal{T} u^{k}-v^{k}-\frac{p^{k}}{\delta}\bigg\|^2+\sigma\|v^{k}\|_{\varphi}\nonumber.
\end{align}
\normalsize
Thus, we have
\footnotesize
\begin{align}
F(w^k)&\geq \bigg(\frac{1}{2}-\frac{\theta^2\|K\|_2^4}{2\delta}\bigg)\|Ku^{k}-f\|^2\nonumber\\
 &\qquad+\frac{\delta}{2}\bigg\|\mathcal{T} u^{k}-v^{k}-\frac{p^{k}}{\delta}\bigg\|^2+\sigma\|v^{k}\|_{\varphi}\nonumber.
\end{align}
\normalsize
From Lemma \ref{descend1}, $\{F(w^k)\}_{k=0,1,2,\ldots}$ is bounded.
This {means}, using the definition of $F(w^k)$, the boundedness of sequences $\{\|v^k\|_{\varphi}\}_{k=0,1,2,\ldots}$,  $\{\mathcal{T} u^{k}-v^{k}-\frac{p^{k}}{\delta}\}_{k=0,1,2,\ldots}$ and $\{Ku^k-f\}_{k=0,1,2,\ldots}$. Taking into account the coercivity of $\varphi$, $\{v^k\}_{k=0,1,2,\ldots}$ is bounded.  Now, taking into account (\ref{boundt1}), $\{p^k\}_{k=0,1,2,\ldots}$ is bounded. Then, $\bigg\{\left(
     \begin{array}{c}
       \mathcal{T} \\
       K \\
     \end{array}
   \right)
u^k\bigg\}_{k=0,1,2,\ldots}$ is bounded, and using the reversibility of $\left(
     \begin{array}{c}
       \mathcal{T} \\
       K \\
     \end{array}
   \right)$, we obtain that $\{u^{k}\}_{k=0,1,2,\ldots}$ is bounded. Therefore, all the components of
the sequence $\{w^k\}_{k=0,1,2,\ldots}$ are bounded, and so the result is proved.
%%%%%%%%%%%%%%%%%%%%%%%%%%%%%%%%%%%%%%%%%%%%%%%%%%%%%%
\section*{Appendix E: Proof of Lemma \ref{relative}}
For the $v^{k+1}$ component, we  have
\footnotesize
\begin{equation}\label{ypart1}
    -\widehat{p}^{k}+\delta (\mathcal{T} u^{k}-v^{k+1})\in\partial \sigma\|v^{k+1}\|_{\varphi}.
\end{equation}
\normalsize
Easy computations on the new defined term $s^{k+1}_v$ give
\footnotesize
\begin{equation}\label{ypart}
    s^{k+1}_v:=(1+\delta)(p^{k+1}-\widehat{p}^{k}) \in \partial_{v} F(w^{k+1}).
\end{equation}
\normalsize
Direct calculations hold
\footnotesize
\begin{align}
    \|s^{k+1}_v\|&\leq(1+\delta)\|p^{k+1}-p^{k}\|+(1+\delta)\|p^{k}-\widehat{p}^{k}\|\nonumber\\
    &=(1+\delta)\|p^{k+1}-p^{k}\|+\alpha(1+\delta)\|p^{k}-p^{k-1}\|\nonumber\\
    &\leq(1+\delta)\theta\|K\|_2^2\|u^{k+1}-u^{k}\|  + \alpha(1+\delta)\theta\|K\|_2^2\|u^{k}-u^{k-1}\|\nonumber.
\end{align}
\normalsize
Thus, we obtain that
\footnotesize
\begin{equation}
\label{eqqv}
    \|s^{k+1}_v\|\leq \gamma_v(\|u^{k+1}-u^{k}\|+\|u^{k}-u^{k-1}\|),
\end{equation}
\normalsize
where the parameter $\gamma_v$ is given by
\footnotesize
$$\gamma_v=\max\bigg\{(1+\delta)\theta\|K\|_2^2, \,
\alpha(1+\delta)\theta\|K\|_2^2\bigg\}.$$
\normalsize

For the $u^{k+1}$ component, we  have
\footnotesize
\begin{equation}\label{xpart1}
  \mathcal{T}^{*}\widehat{p}^{k}-\delta \mathcal{T}^{*}(\mathcal{T} u^{k+1}-v^{k+1})=K^*(Ku^{k+1}-f).
\end{equation}
\normalsize
With (\ref{xpart1}), we have the new term $s^{k+1}_u$
\footnotesize
 \begin{align}\label{xpart}
    s^{k+1}_u&:=\mathcal{T}^{*}\widehat{p}^{k}-\mathcal{T}^{*}p^{k+1}+
    \frac{7\alpha^2\theta^2\|K\|_2^4}{2\delta}(u^{k+1}-u^k)\in \partial_{u} F(w^{k+1}).
\end{align}
\normalsize
It is easy to obtain
\footnotesize
\begin{eqnarray}
\label{eqqu}
    \|s^{k+1}_u\|\leq \gamma_u(\|u^{k+1}-u^{k}\|+\|u^{k}-u^{k-1}\|),
\end{eqnarray}
\normalsize
where \footnotesize
$$\gamma_u=\max\bigg\{\theta\|\mathcal{T}\|_2\|K\|_2^2+\frac{7\alpha^2\theta^2\|K\|_2^4}{2\delta}, \,
\alpha\theta\|\mathcal{T}\|_2\|K\|_2^2\bigg\}.$$
\normalsize
Obviously, it holds
\footnotesize
\begin{equation}
\label{eqqp}
    s^{k+1}_p:=-\mathcal{T} u^{k+1}+v^{k+1}\in \partial_{p} F(w^{k+1}).
\end{equation}
\normalsize
From the scheme, we can easily see
\footnotesize
\begin{align}
    \|s^{k+1}_p\|&=\frac{1}{\delta}\|p^{k+1}-\widehat{p}^k\|\leq\frac{\|p^{k+1}-p^k\|}{\delta}+\frac{\alpha\|p^{k}-p^{k-1}\|}{\delta}\nonumber\\
    &\leq\gamma_p(\|u^{k+1}-u^{k}\|+
    \|u^{k}-u^{k-1}\|), \nonumber
\end{align}
\normalsize
where $\gamma_p=\max\bigg\{\frac{\theta\|K\|_2^2}{\delta}
, \,
\frac{\alpha\theta\|K\|_2^2}{\delta}\bigg\}.$
Therefore, we have
\footnotesize
\begin{equation}
\label{eqqx}
    \|\nabla_{x}F(w^{k+1})\|\leq \frac{7\alpha^2\theta^2\|K\|_2^4}{\delta} (\|d^{k+1}-d^{k}\|+\|u^{k}-u^{k-1}\|).
\end{equation}
\normalsize

Letting $s^{k+1}:=(s^{k+1}_v,s^{k+1}_u,s^{k+1}_p,\nabla_{x}F(w^{k+1}))$, we have $s^{k+1}\in \partial F(w^{k+1})$. With all the  above equations (\ref{eqqv}), (\ref{eqqu}), (\ref{eqqp}), and (\ref{eqqx}), we obtain the global bound
\footnotesize
\begin{align}
    \|s^{k+1}\|&\leq(\gamma_u+\gamma_v+\gamma_p+\frac{7\alpha^2\theta^2\|K\|_2^4}{\delta})\nonumber\\
    & \qquad \times(\|u^{k+1}-u^{k}\|+\|u^{k}-u^{k-1}\|).
\end{align}
\normalsize
Denoting
$\gamma=\gamma_u+\gamma_v+\gamma_p+\frac{7\alpha^2\theta^2\|K\|_2^4}{\delta},$
and then, we finish the proof.

%%%%%%%%%%%%%%%%%%%%%%%%%%%%%%%%%%%%%%%%
\section*{Appendix F: Proof of Theorem \ref{th1}}
For any cluster point $(u^*,v^*,p^*)$, there exists $\{k_j\}_{j=0,1,2,\ldots}$ such that $\lim_{j}(u^{k_j},v^{k_j},p^{k_j})=(u^*,v^*,p^*)$. Then, from Lemma \ref{descend1}, we have
$
\lim_{j}(u^{k_j+1},v^{k_j+1},p^{k_j+1})=\lim_{j}(\widehat{u}^{k_j},\widehat{v}^{k_j},\widehat{p}^{k_j})=(u^*,v^*,p^*).
$
From the scheme of Algorithm~\ref{algo}, we have the following conditions
\footnotesize
\begin{align}
 -\widehat{p}^{k_j}+\delta(\mathcal{T} u^{k_j}-v^{k_j+1})&\in \partial \sigma \|v^{k_j+1}\|_{\varphi},\nonumber\\
\mathcal{T}^*\widehat{p}^{k_j}-\delta \mathcal{T}^{*}(\mathcal{T} u^{k_j+1}-v^{k_j+1})&=K^*(Ku^{k_j+1}-f),\nonumber\\
 p^{k_j+1}&=\widehat{p}^{k_j}-\delta(\mathcal{T} u^{k_j+1}-\widehat{v}^{k_j+1}).\nonumber
\end{align}
\normalsize
Letting $j\rightarrow+\infty$, with Proposition \ref{sublimit}, we are then led to
\footnotesize
\begin{eqnarray}
 -p^{*}&\in& \partial \sigma \|v^{*}\|_{\varphi},\nonumber\\
\mathcal{T}^*p^{*}&=&K^*(Ku^{*}-f),\nonumber\\
\mathcal{T} u^{*}&=&v^*.\nonumber
\end{eqnarray}
\normalsize
Finally, from Proposition~\ref{criL}, we obtain that $(u^*,v^*,p^*)$ is a critical point of $\mathcal{L}^{\varphi}_{\delta}$.

%%%%%%%%%%%%%%%%%%%%%%%%%%%%%%%%%%%%%%%%%%%%%%%%

\section*{Appendix G: Proof of Lemma \ref{station}}
(1) From Lemma~\ref{bound} the sequences $\{(u^k,v^k,p^k)\}_{k=0,1,2,\ldots}$  and $\{w^{k}\}_{k=0,1,2,\ldots}$ are bounded. Thus, $\mathcal{M}$ is nonempty. Assume that $w^*\in \mathcal{M}$, then from the definition there exists
a subsequence $w^{k_i}=(u^{k_i},v^{k_i},p^{k_i},u^{k_i-1})\rightarrow w^*=(u^{*},v^{*},p^{*},u^{\bullet})$. In that case, from Lemma~\ref{descend1}, we also have $u^{k_i-1}\rightarrow u^*$. That is also $u^{\bullet}=u^*$, i.e., $w^*=(u^{*},v^{*},p^{*},u^{*})$.   Besides, from Lemma~\ref{relative} and Lemma~\ref{descend1}, we have $s^{k_i}\in \partial F(w^{k_i})$ and
$s^{k_i}\rightarrow \mathbf{0}$.  Finally, Proposition~\ref{sublimit} indicates that $\mathbf{0}\in \partial F(w^*)$, i.e. $w^*\in \textrm{crit} (F)$.

(2) This item follows as a consequence of the definition of the limit point set (Definition~\ref{deflim}).

(3) In Eq.~(\ref{s1temp2}), by taking limits on $i$, we obtain
$
    \lim\sup_{i}\|v^{k_{i}}\|_{\varphi}\leq\|v^{*}\|_{\varphi}.
$
And with the closedness of $\varphi$, we have the other bound
   $ \lim\inf_{i}\|v^{k_{i}}\|_{\varphi}\geq\|v^{*}\|_{\varphi}$.
That means
$\lim_{i}\|v^{k_{i}}\|_{\varphi}=\|v^{*}\|_{\varphi}.
$
Using the continuity of the rest of functions that defines $F$, we have
$
    \lim_{i}F(w^{k_i})=F(w^*).
$
 On the other hand, Lemma~\ref{descend1} indicates $\{F(w^k)\}_{k=0,1,2,\ldots}$ is decreasing. Noting the boundedness of $\{F(w^k)\}_{k=0,1,2,\ldots}$, this sequence has a lower bound. Thus, $\{F(w^k)\}_{k=0,1,2,\ldots}$ is convergent. And then, we have
$
    \lim_{k}F(w^{k})=\lim_{i}F(w^{k_i})=F(w^*).
$

\section*{Appendix H: Proof of Theorem \ref{th2}}
%As $\varphi$ is a semi-algebraic function, then the auxiliary function $F$ (\ref{fun}) is also semi-algebraic.
From Lemma~\ref{station}, $F$ is constant on $\mathcal{M}$.
Let $w^*=(u^*,v^*,p^*,x^*)$ be a stationary point of $\{w^k\}_{k=0,1,2,\ldots}$. Then, from the definition of $w^k$ we have
$
    u^*=x^*.
$
Also from Lemma~\ref{station}, we have $\textrm{dist}(w^k,\mathcal{M})<\varepsilon$ and $F(w^k)<F(w^*)+\eta$ for any $k>K$, for some $K$. Hence, from Lemma~\ref{con}, we have
\footnotesize
\begin{equation}
    \textrm{dist}(\textbf{0},\partial F(w^k))\rho'(F(w^k)-F(w^*))\geq 1,
\end{equation}
\normalsize
that together with Lemma~\ref{relative} give us
\footnotesize
\begin{align}
    &\frac{1}{\rho'(F(w^k)-F(w^*))}\leq\textrm{dist}(\textbf{0},\partial F(w^k))\nonumber\\
    &\quad\quad\leq \gamma(\|u^{k+1}-u^{k}\|+\|u^{k}-u^{k-1}\|).
\end{align}
\normalsize
Then, the concavity of $\rho$ yields that
\footnotesize
\begin{align}
    &F(w^k)-F(w^{k+1})= F(w^k)-F(w^*)-[F(w^{k+1})-F(w^*)]\nonumber\\
    &\quad\leq\frac{\rho[F(w^k)-F(w^*)]-\rho[F(w^{k+1})-F(w^*)]}{\rho'[F(w^k)-F(w^*)]}\nonumber\\
    &\quad \leq\gamma(\|u^{k+1}-u^{k}\|+\|u^{k}-u^{k-1}\|)\nonumber\\
    &  \quad \times \{\rho[F(w^k)-F(w^*)]-\rho[F(w^{k+1})-F(w^*)]\}.
\end{align}
\normalsize
Using Lemma~\ref{relative}, we have
\footnotesize
\begin{align}
   &\widehat{h}\|u^{k+1}-u^k\|^2\leq \gamma(\|u^{k+1}-u^{k}\|+\|u^{k}-u^{k-1}\|) \nonumber\\
   &\quad \times \{\rho[F(w^k)-F(w^*)]-\rho[F(w^{k+1})-F(w^*)]\},
\end{align}
\normalsize
which is equivalent to
\footnotesize
\begin{align}
   &\frac{\widehat{h}}{\gamma}\|u^{k+1}-u^k\|\leq2\cdot\frac{1}{2}\sqrt{\frac{\widehat{h}}{\gamma}}\sqrt{\|u^{k+1}-u^{k}\|+\|u^{k}-u^{k-1}\|}\nonumber\\
   &\quad\times\sqrt{\rho[F(w^k)-F(w^*)]-\rho[F(w^{k+1})-F(w^*)]}.
\end{align}
\normalsize
Using the Schwartz's inequality, we then derive
\footnotesize
\begin{align}\label{core}
    &\frac{\widehat{h}}{\gamma}\|u^{k+1}-u^k\|\leq\frac{\widehat{h}}{4\gamma}(\|u^{k+1}-u^{k}\|+\|u^{k}-u^{k-1}\|) \nonumber\\
    &\quad +\{\rho[F(w^k)-F(w^*)]-\rho[F(w^{k+1})-F(w^*)]\}.
\end{align}
\normalsize
Summing (\ref{core}) from $K$ to $K+j$ yields that
\footnotesize
$$
\begin{array}{l}
   \displaystyle\frac{\widehat{h}}{2\gamma} \displaystyle\sum_{k=K}^{K+j-2}\|u^{k+1}-u^k\|+\displaystyle\frac{3\widehat{h}}{4\gamma}\|u^{K+j+1}-u^{K+j}\| \\[9.pt]
    \qquad \leq\rho[F(w^K)-F(w^*)]
   -\varphi[F(w^{K+j+1})-F(w^*)]\\[5.pt]
   \qquad  \qquad +\displaystyle\frac{\widehat{h}}{4\gamma}\|u^{K-1}-u^{K-2}\|.
\end{array}
$$
\normalsize
Letting $j\rightarrow+\infty$, and applying Lemma~\ref{descend1}, we have
$
   \frac{\widehat{h}}{2\gamma} \sum_{k=K}^{+\infty}\|u^{k+1}-u^k\|<+\infty.
$
Then, $\{u^k\}_{k=0,1,2,\ldots}$ is convergent. And noting that $u^*$ is a stationary point of $\{u^k\}_{k=0,1,2,\ldots}$ we also have that $\{u^k\}_{k=0,1,2,\ldots}$  converges to $u^*$. With (\ref{condition1}), the convergence of $\{u^k\}_{k=0,1,2,\ldots}$ indicates the convergence of $\{p^k\}_{k=0,1,2,\ldots}$. And then the sequence $\{\hat{p}^k\}_{k=0,1,2,\ldots}$ is also convergent. With the identity $p^{k+1}=\widehat{p}^k-\delta(\mathcal{T} u^{k+1}-v^{k+1})$, $\{v^k\}_{k=0,1,2,\ldots}$ is convergent.
That means   $(u^k,v^k,p^k)$ converges to some point $(u^*,v^*,p^*)$. Finally, from Theorem~\ref{th1}, we have that $(u^*,v^*,p^*)$ is a critical point of $\mathcal{L}^{\varphi}_{\delta}$.

\ifCLASSOPTIONcaptionsoff
  \newpage
\fi

\end{document}